\documentclass[12pt]{amsart}
\usepackage{amsthm, amssymb}
\usepackage{amsfonts, amscd}
\usepackage{epsfig,multicol}
\theoremstyle{plain} \numberwithin{equation}{section}
\newtheorem{theo}{Theorem}[section]
\newtheorem{coro}[theo]{Corollary}
\newtheorem{prop}[theo]{Proposition}
\newtheorem{lemm}[theo]{Lemma}
\newtheorem*{prob}{Problem}
\newtheorem*{theo*}{Theorem}

\theoremstyle{definition}
\newtheorem*{defi}{Definition}
\newtheorem{exam}[theo]{Example}
\newtheorem*{rema}{Remark}
\newtheorem*{conj}{Conjecture}

\def\Z{\mathbb Z}
\def\C{\mathbb C}
\def\R{\mathbb R}

\def\a{\mathbf a}
\def\b{\mathbf b}

\def\e{\mathbf e}

\def\mZ{\mathcal{Z}}

\DeclareMathOperator{\Hom}{Hom}

\begin{document}
\title{Buchstaber invariants of skeleta of a simplex}
\author[Y. Fukukawa]{Yukiko Fukukawa}
\author[M. Masuda]{Mikiya Masuda}
\address{Department of Mathematics, Osaka City
University, Sumiyoshi-ku, Osaka 558-8585, Japan.}
\email{m09sa024@ex.media.osaka-cu.ac.jp}
\email{masuda@sci.osaka-cu.ac.jp}

\date{\today}
\thanks{The second author was partially supported by Grant-in-Aid for Scientific Research 19204007}

\begin{abstract}
A moment-angle complex $\mZ_K$ is a compact topological space associated with a finite simplicial complex $K$.  
It is realized as a subspace of a polydisk $(D^2)^m$, where $m$ is the number of vertices in $K$ and $D^2$ is the unit disk 
of the complex numbers $\C$, and the natural 
action of a torus $(S^1)^m$ on $(D^2)^m$ leaves $\mZ_K$ invariant.  The Buchstaber invariant $s(K)$ of $K$ 
is the maximum integer for which there is a subtorus of rank $s(K)$ acting on $\mZ_K$ freely.  

The story above goes over the real numbers $\R$ in place of $\C$ and a real analogue of the Buchstaber invariant, 
denoted $s_\R(K)$, can be defined for $K$ and $s(K)\leqq s_\R(K)$.  
In this paper we will make some computations of $s_\R(K)$ when $K$ is a skeleton of a simplex.  
We take two approaches to find $s_\R(K)$ and the latter one turns out to be a problem of integer linear programming 
and of independent interest.   
\end{abstract}

\maketitle

\section{Introduction}
Davis and Januszkiewicz (\cite{da-ja91}) initiated the study of topological analogue of toric geometry and 
introduced a compact topological space $\mZ_K$ associated with a finite simplicial complex $K$.  
Then Buchstaber and Panov (\cite{bu-pa02})  intensively studied the topology of $\mZ_K$ by realizing it in a polydisk $(D^2)^m$,   
where $m$ is the number of vertices in $K$ and $D^2$ is the unit disk of the complex numbers $\C$, 
and noted that $\mZ_K$ is a deformation retract of the complement  
of the union of coordinate subspaces in $\C^m$ associated with $K$.  
They named $\mZ_K$ a \emph{moment-angle complex} associated with $K$.  Although the construction of $\mZ_K$ is simple, 
the topology of $\mZ_K$ is complicated in general and the space $\mZ_K$ is getting more attention of topologists, see 
\cite{cont08}.    

The coordinatewise multiplication of a torus $(S^1)^m$ on $\C^m$, where $S^1$ is the unit circle of $\C$, leaves $\mZ_K$ 
invariant.  
The action of $(S^1)^m$ on $\mZ_K$ is not free but its restriction to a certain 
subtorus of $(S^1)^m$ can be free.  The maximum integer $s(K)$ for which there is a subtorus of dimension $s(K)$ 
acting freely on $\mZ_K$ is a combinatorial invariant and called the \emph{Buchstaber invariant} of $K$.  
When $K$ is of dimension $n-1$, $s(K)\leqq m-n$ and Buchstaber (\cite{buch08}, \cite{bu-pa02}) asked 

\begin{prob}
Find a combinatorial description of $s(K)$. 
\end{prob}

If $P$ is a simple convex polytope of dimension $n$, then its dual $P^*$ is a simplicial polytope and the boundary 
$\partial P^*$ of $P^*$ is a simplicial complex of dimension $n-1$.  The Buchstaber invariant $s(P)$ of $P$ is then defined to be 
$s(\partial P^*)$.  We note that $s(P)=m-n$, where $m$ is the number of vertices of $P^*$, if and only if there is a quasitoric 
manifold over $P$.  
We refer the reader to \cite{aize09} and \cite{erok09} for some properties and computations on $s(P)$ and $s(K)$.  
The reader can also find some results on them in \cite[Theorem 6.6]{buch08}.  
 
The story mentioned above goes over the real numbers $\R$ in place of $\C$.  In this case, the moment-angle complex 
$\mZ_K$ is replaced by a \emph{real moment-angle complex} $\R\mZ_K$ and the torus $(S^1)^m$ is replaced by 
a 2-torus $(S^0)^m$ 
where $S^0=\{\pm 1\}$.  Then a real analogue of the Buchstaber invariant can be defined for $K$, which we denote by $s_\R(K)$. 
Namely $s_\R(K)$ is the maximum integer for which there is a 2-subtorus of rank $s_\R(K)$ 
acting freely on $\R\mZ_K$.  The complex conjugation on $\C$ induces an involution on $\mZ_K$ with $\R\mZ_K$ as 
the fixed point set and this implies that $s(K)\leqq s_\R(K)$.  
 
In this paper we make some computations of  
$s_\R(K)$ when $K$ is a skeleton of a simplex.   Let $\Delta^{m-1}_r$ be the $r$-skeleton of the $(m-1)$-simplex.  Then it follows from 
the definition of $\R\mZ_K$ (see \cite[p.98]{bu-pa02}) that 
\begin{equation} \label{eqZ}
\R\mZ_{\Delta^{m-1}_{m-p-1}}=\bigcup (D^1)^{m-p}\times (S^0)^{p} \subset (D^1)^m
\end{equation}
where $D^1$ is the interval $[-1,1]$ in $\R$ so that $S^0$ is the boundary of $D^1$ and 
the union is taken over all $m-p$ products of $D^1$ in $(D^1)^m$.  
We denote the invariant $s_\R(\Delta^{m-1}_{m-p-1})$ simply by $s_\R(m,p)$.  
The moment-angle complex $\R\mZ_{\Delta^{m-1}_{m-p-1}}$ is sitting in the complement 
$U_\R(m,p)$ of the union of all coordinate subspaces of dimension $p-1$ in $\R^m$ and 
$s_\R(m,p)$ may be thought of as the maximal integer for which there is a 2-subtorus of rank $s_\R(m,p)$ 
acting freely on $U_\R(m,p)$.

We easily see $s_\R(m,0)=0$ and assume $p\geqq 1$.  
We take two approaches to find $s_\R(m,p)$ and here is a summary of the results obtained from the first approach developed in 
Section~\ref{sect:3}. 

\begin{theo*}
Let $1\leqq p\leqq m$.  
\begin{enumerate}
\item $1\leqq s_\R(m,p)\leqq p$ and $s_\R(m,p)=p$ if and only if $p=1,m-1,m$. 
\item $s_\R(m,p)$ increases as $p$ increases but decreases as $m$ increases. 
\item If $m-p$ is even, then $s_\R(m,p)=s_\R(m+1,p)$.
\item $s_\R(m+1,m-2)=s_\R(m,m-2)=[m-\log_2(m+1)]$ for $m\geqq 3$, where $[r]$ for a real number $r$ denotes the greatest integer less 
than or equal to $r$. 
\end{enumerate}
\end{theo*}

\begin{rema}
It is easy to prove (1) and (2) above.  After we finished writing the first version of this paper, we learned from N. Erokhovets 
that (4) was also obtained by A. Aizenberg \cite{aize09}, see also \cite{erok09}.  
\end{rema}

It seems difficult to find a computable description of $s_\R(m,p)$ in terms of $m$ and $p$ in general.  From Section~\ref{sect:4} 
we take another approach to find $s_\R(m,p)$, that is, 
we investigate values of $m$ and $p$ for which $s_\R(m,p)$ is a given positive integer $k$.   
It turns out that $s_\R(m,p)=1$ if and only if $m\geqq 3p-2$ (Theorem~\ref{theo3}) and that  
there is a non-negative integer $m_k(b)$ associated to integers $k\geqq 2$ and $b\geqq 0$ such that 
\[
\begin{split}
s_\R(m,p)&\geqq k \text{ if and only if $m\leqq m_k(p-1)$,}\\
\text{in other words,}&\text{ since $s_\R(m,p)$ decreases as $m$ increases,}\\
s_\R(m,p)&=k \text{ if and only if $m_{k+1}(p-1)<m\leqq m_k(p-1)$.}  
\end{split}
\]
Therefore, finding $s_\R(m,p)$ is equivalent to finding $m_k(p-1)$ for all $k$.  
In fact, $m_k(b)$ is the maximum integer which the linear function $\sum_{v\in (\Z/2)^k\backslash\{0\}} a_v$ 
takes on lattice points $(a_v)$ in $\R^{2^k-1}$ satisfying these $(2^k-1)$ inequalities 
\[
\sum_{(u,v)=0} a_v\leqq  b\quad \text{for each $u\in (\Z/2)^k\backslash\{0\}$}
\]
and $a_v\geqq 0$ for every $v$, where $\Z/2=\{0,1\}$ and $(\ ,\ )$ denotes the standard scalar product on $(\Z/2)^k$. 
Finding $m_k(b)$ is a problem of integer linear programming and of independent interest.  
Here is one of the main results on $m_k(b)$.  

\begin{theo*}[Theorem~\ref{theo:main}] 
Let $b=(2^{k-1}-1)Q+R$ with non-negative integers $Q,R$ with $0\leqq R\leqq 2^{k-1}-2$.  
We may assume that $2^{k-1}-2^{k-1-\ell}\leqq R\leqq 2^{k-1}-2^{k-1-(\ell+1)}$ for some $0\leqq \ell\leqq k-2$. 
Then   
\[
(2^k-1)Q+R+2^{k-1}-2^{k-1-\ell}\leqq m_k(b)\leqq (2^k-1)Q+2R,
\]
and 
the lower bound is attained if and only if $R-(2^{k-1}-2^{k-1-\ell})\leqq k-\ell-2$ and 
the upper bound is attained if and only if $R=2^{k-1}-2^{k-1-\ell}$.  
\end{theo*}
More explicit values of $m_k(b)$ can be found in Sections~\ref{sect:6} and~\ref{sect:56}.  
In particular, $m_k(b)$ is completely determined for $k=2,3,4$, see Example~\ref{ex:k=34}, so that 
one can find for which values of $m$ and $p$ we have $s_\R(m,p)\geqq k$ for $k=2,3,4$.  
The equivalent results are obtained in \cite{erok09} for $k=2,3$.   

All of our computations support a conjecture that 
$$m_k((2^{k-1}-1)Q+R)=(2^k-1)Q+m_k(R)$$
would hold for any $Q$ and $R$.  
This is equivalent to $m_k(b+2^{k-1}-1)=m_k(b)+2^k-1$ for any $b$ and  we prove in Section~\ref{sect:8} that 
the latter identity holds when $b$ is large.   
 

The authors thank Suyoung Choi for his help to find $k\times m$ matrices which realize 
$s_\R(m,p)=k$ for small values of $m$ and $p$.   They also thank Nickolai Erokhovets
for helpful comments on an earlier version of the paper.

\section{Some properties and computations of $s_\R(m,p)$} \label{sect:3}

In this section we translate our problem to a problem of linear algebra, deduce some properties of 
$s_\R(m,p)$ and make some computations of $s_\R(m,p)$.  

The real moment-angle complex $\R\mZ_{\Delta^{m-1}_{m-p-1}}$ in \eqref{eqZ} with $p=0$ 
is the disk $(D^1)^m$. Since the action of $(S^0)^m$ on $(D^1)^m$ 
has a fixed point, that is the origin, we have 
\begin{equation} \label{sm0}
s_\R(m,0)=0.
\end{equation}
Another extreme case is when $p=m$.  Since $\R\mZ_{\Delta^{m-1}_{m-p-1}}$ in \eqref{eqZ} with $p=m$ is 
$(S^0)^m$,  we have 
\begin{equation} \label{smm}
s_\R(m,m)=m.
\end{equation}
In the following we assume $p\geqq 1$.  

\begin{lemm} \label{lemm1}
Let $A=(\a_1,\dots,\a_m)$ be a $k\times m$ matrix with entries in $\Z/2$ 
and let $\rho_A\colon (S^0)^k\to (S^0)^m$ be a homomorphism defined by $\rho_A(g)=(g^{\a_1},\dots,g^{\a_m})$, where 
$g^\a=\prod_{i=1}^k g_i^{a^i}$ for $g=(g_1,\dots,g_k)\in (S^0)^k$ and a column vector 
$\a=(a^1,\dots,a^k)^T$ in $(\Z/2)^k$.  
Then the action of $(S^0)^k$ on $\R\mZ_{\Delta^{m-1}_{m-p-1}}$ in \eqref{eqZ} 
through $\rho$ is free if and only if any $p$ column vectors in $A$ span $(\Z/2)^{k}$.
\end{lemm}

\begin{proof}
The action of $(S^0)^k$ on $\R\mZ_{\Delta^{m-1}_{m-p-1}}$ 
through $\rho_A$ leaves each subspace $(D^1)^{m-p}\times (S^0)^{p}$ in \eqref{eqZ} invariant and 
the action on $\R\mZ_{\Delta^{m-1}_{m-p-1}}$ 
 is free if and only if it is free on each $(D^1)^{m-p}\times (S^0)^{p}$.  
The latter is equivalent to the action being free on each $\{0\}\times (S^0)^{p}$ and this is equivalent to 
$\rho$ composed with the projection from $(S^0)^m$ onto $(S^0)^{p}$ being injective.  
This is further equivalent to a matrix formed from any $p$ column vectors in 
$A$ being of full rank (that is $k$), which is equivalent to the last statement in the lemma. 
\end{proof}

Since any rank $k$ subgroup of $(S^0)^m$ is obtained as $\rho_A((S^0)^k)$ for some $A$ in Lemma~\ref{lemm1}, 
Lemma~\ref{lemm1} implies  

\begin{coro} \label{coro1}
The invariant $s_\R(m,p)$ is the maximum integer $k$ for which 
there exists a $k\times m$ matrix $A$ with entries in $\Z/2$ such that any $p$ column vectors in $A$ span $(\Z/2)^k$.  
\end{coro}

Here are some properties of $s_\R(m,p)$. 

\begin{prop} \label{prop2}
\begin{enumerate}
\item $1\leqq s_\R(m,p)\leqq p$ for $p\geqq 1$.  In particular, $s_\R(m,1)=1$. 
\item $s_\R(m,p)\leqq s_\R(m,p')$ if $p\leqq p'$. 
\item $s_\R(m,p)\geqq s_\R(m',p)$ if $m\leqq m'$.
\end{enumerate}
\end{prop}

\begin{proof}
The inequality 
(1) is obvious from Corollary~\ref{coro1} and the inequality (2) follows from the fact that 
if $p'\geqq p$, then $\R\mZ_{\Delta^{m-1}_{m-p-1}}$ in \eqref{eqZ} contains 
$\R\mZ_{\Delta^{m-1}_{m-p'-1}}$ as an invariant subspace. 

Let $m'\geqq m$ and set $k=s_\R(m',p)$.  Then there is a $k\times m'$ matrix $A'$ with entries $\Z/2$ such that any $p$ 
column vectors in $A'$ span $(\Z/2)^k$.  Let $A$ be a $k\times m$ matrix formed from arbitrary $m$ column vectors in $A'$. 
Since any $p$ column vectors in $A$ span $(\Z/2)^k$, it follows from Corollary~\ref{coro1} that $s_\R(m,p)\geqq k=s_\R(m',p)$.   
 \end{proof}


We denote by $\{\e_1,\dots,\e_k\}$ the standard basis of $(\Z/2)^k$.  

\begin{theo} \label{theo1}
$s_\R(m,m-1)=m-1$ for $m\geqq 2$.
\end{theo}

\begin{proof}
We have $s_\R(m,m-1)\leqq m-1$ by Proposition~\ref{prop2} (1). 
On the other hand, any $m-1$ column vectors in an $(m-1)\times m$ matrix $A=(\e_1,\dots,\e_{m-1}, \sum_{i=1}^{m-1}\e_i)$ 
span $(\Z/2)^{m-1}$, 
so $s_\R(m,m-1)\geqq m-1$ by Lemma~\ref{lemm1}.  
\end{proof}

If $A$ is a $k\times m$ matrix with entries in $\Z/2$ which realizes $s_\R(m,p)=k$, then $A$ must be of full rank (that is $k$); so 
we may assume that the first $k$ column vectors in $A$ are 
linearly independent if necessary by permuting columns and moreover that they  
are $\e_1,\dots,\e_k$ by multiplying $A$ by an invertible matrix of size $k$ from the left.  

\begin{lemm} \label{lemm2}
$s_\R(m,p)\leqq p-1$ when $2\leqq p\leqq m-2$.
\end{lemm}

\begin{proof}
Since $s_\R(m,p)\leqq p$ by Proposition~\ref{prop2} (1), it suffices to prove that $s_\R(m,p)\not=p$ when $2\leqq p\leqq m-2$. 
Suppose $s_\R(m,p)=p$ and let $A$ be a $p\times m$ matrix $(\e_1,\dots,\e_{p}, \a_{p+1},\dots,\a_m)$ 
which realizes $s_\R(m,p)=p$.  Then all $\a_j$'s for $j=p+1,\dots,m$ must be equal to $\sum_{i=1}^{p} \e_i$ because 
any $p-1$ vectors from $\e_1,\dots,\e_{p}$ together with one $\a_j$ span $(\Z/2)^p$. 
The number of $\a_j$'s is more than one as $p\leqq m-2$, so $p$ column vectors in $A$ containing more than 
one $\a_j$ do not span $(\Z/2)^p$, which is a contradiction.  
\end{proof}

\begin{theo} \label{mm+1}
If $m-p$ is even, then $s_\R(m,p)=s_\R(m+1,p)$.
\end{theo}

\begin{proof}
The original proof of this theorem was rather long.  Below is a much simpler proof due to 
Nickolai Erokhovets.  We thank him for sharing his argument.  

Since $s_\R(m,0)=0$ for any $m$ by \eqref{sm0}, we may assume $p\geqq 1$ so that 
we can use Corollary~\ref{coro1}.  
Suppose that $m-p$ is even and set $s_\R(m,p)=k$. 
Since $s_\R(m,p)$ decreases as $m$ increases by Proposition~\ref{prop2} (3), it suffices to show that 
there is a $k\times (m+1)$ matrix in which any $p$ column vectors span $(\Z/2)^k$.  

Let $A=(\a_1,\dots,\a_{m})$ be a $k\times m$ matrix  which realizes $s_\R(m,p)=k$. 
Set $\b=\sum_{i=1}^m\a_i$ and consider a $k\times(m+1)$ matrix $B=(\a_1,\dots,\a_m,\b)$. 
We shall prove that any $p$ column vectors in $B$ span $(\Z/2)^k$.  
If $\b$ is not a member of the $p$ column vectors, then all of them are in $A$ so that they span $(\Z/2)^k$ by the choice of $A$.  
Therefore we may assume that $\b$ is a member of the $p$ column vectors.  
If the $p-1$ column vectors except $\b$, say $\a_{i_1},\dots,\a_{i_{p-1}}$, span $(\Z/2)^k$, then we have nothing to do. 
Suppose that the $p-1$ column vectors do not span $(\Z/2)^k$.  Then they span a codimension 
1 subspace, say $V$, of $(\Z/2)^k$ because $\a_{i_1},\dots,\a_{i_{p-1}}$ are in $A$ and any $p$ column vectors 
in $A$ span $(\Z/2)^k$ by the choice of $A$.  This shows that if $f$ is a homomorphism from $(\Z/2)^k$ to $\Z/2$  
whose kernel is $V$, then $f(\a_{i_j})=0$ for $j=1,\dots,p-1$ and $f(\a_\ell)=1$ for any $\ell$ different from 
$i_1,\dots,i_{p-1}$.  It follows that 
\[
f(\b)=f(\sum_{i=1}^m\a_i)=m-(p-1)=1 \in \Z/2
\]
where we used the assumption on $m-p$ being even at the last identity. Therefore $\b$ is not contained in $V$ 
so that the $p$ column vectors $\a_{i_1},\dots,\a_{i_{p-1}},\b$ span $(\Z/2)^k$.  
This completes the proof of the theorem.
\end{proof}

If we take $p=m-2\geqq 4$ in Lemma~\ref{lemm2}, we have $s_\R(m,m-2)\leqq m-3$ for $m\geqq 4$.  In fact, 
$s_\R(m,m-2)$ is given as follows.  

\begin{theo} \label{theo2}
$s_\R(m+1,m-2)=s_\R(m,m-2)=[ m-\log_2(m+1)]$ for $m\geqq 3$. 
\end{theo}

\begin{proof}
The first identity follows from Theorem~\ref{mm+1}, so it suffices to prove the second identity.  

Set $s_\R(m,m-2)=k$ and let $A=(\e_1,\dots,\e_k,\a_{k+1},\dots,\a_m)$ be a matrix which realizes $s_\R(m,m-2)=k$.  
Then any $m-2$ column vectors in $A$ span $(\Z/2)^k$.  This means that for each $i=1,\dots,k$ 
the set 
\[
A(i):=\{\ell \mid \text{the $i$-th component of $\a_\ell$ is $1$}\}\subset \{k+1,\dots,m\}
\]
contains at least two elements because if $A(i)$ consists of only one element, say $\ell$, for some $i$, then 
the $m-2$ column vectors in $A$ except $\e_i$ and $\a_\ell$ will not generate a vector with $1$ at the 
$i$-th component.  Another constraint on $A(i)$'s is that they are mutually distinct because if $A(i)=A(j)$ for some $i$ and $j$ in 
$\{1,\dots,k\}$, then $m-2$ column vectors in $A$ except $\e_i$ and $\e_j$ will not generate $\e_i$ 
and $\e_j$.  
Conversely, if $A(i)$ contains at least two elements for each $i$ and $A(i)$'s are mutually distinct, then 
any $m-2$ column vectors in $A$ span $(\Z/2)^k$.  

The number of subsets of $\{k+1,\dots,m\}$ which contain at least two elements is given by 
\[
\sum_{n=2}^{m-k}\binom{m-k}{n}=2^{m-k}-1-m+k.
\]
Since the number of $A(i)$'s is $k$, the argument above shows that $k$ should be the maximum integer which satisfies 
\[
\text{$k\leqq 2^{m-k}-1-m+k$,\quad i.e.,\quad $k\leqq m-\log_2(m+1)$.}
\]  
This proves the theorem.  
\end{proof}

\section{Another approach to compute $s_\R(m,p)$} \label{sect:4}

We know $s_\R(m,p)=p$ when $p=0,1$.  So we will assume $p\geqq 2$ in the following.  
It seems difficult to find a computable description of $s_\R(m,p)$ in terms of $m$ and $p$ in general.  
Hereafter we take a different approach to find values of $s_\R(m,p)$ for $p\geqq 2$, i.e. we find values of $m$ and $p$ 
for which $s_\R(m,p)$ is a given positive integer $k$.  We begin with 

\begin{theo} \label{theo3}
$s_\R(m,p)=1$ if and only if $m\geqq 3p-2$, 
in other words, $s_\R(m,p)\geqq 2$ if and only if $m\leqq 3(p-1)$.  
\end{theo}

\begin{proof}
Since $s_\R(m,p)$ decreases as $m$ increases by Proposition~\ref{prop2} (3), it suffices to show 
\begin{enumerate}
\item $s_\R(3(p-1),p)\geqq 2$, and 
\item $s_\R(3p-2,p) =1$.
\end{enumerate}

\medskip
\noindent  
{Proof of (1)}.   
Let $A$ be a $2\times 3(p-1)$ matrix formed from $p-1$ copies of $(\e_1,\e_2,\e_1+\e_2)$.  
Then any $p$ column vectors in $A$ span $(\Z/2)^2$, 
which means $s_\R(3(p-1),p)\geqq 2$. 

\medskip
\noindent
{Proof of (2)}.  Suppose that $s_\R(3p-2,p)\geqq 2$.  Then there is a $2\times (3p-2)$ matrix $A$ 
such that any $p$ column vectors in $A$ span $(\Z/2)^2$.  Let $\e_i$ (resp. $\e_1+\e_2$) appear $a_i$ 
(resp. $a_{12}$) times in $A$.  Then 
\begin{equation} \label{eq1}
a_1+a_2+a_{12}=3p-2
\end{equation}
 and inequalities
\[
a_i\leqq p-1 \quad\text{for $i=1,2$}\qquad \text{and}\qquad a_{12}\leqq p-1
\]
must be satisfied for any $p$ column vectors in $A$ to span $(\Z/2)^2$. 
These inequalities imply that $a_1+a_2+a_{12}\leqq 3p-3$ which contradicts \eqref{eq1}.
\end{proof}

The above argument can be developed for general values of $k$ with $s_\R(m,p)\geqq k$. 
Let $(\ ,\ )$ be the standard bilinear form on $(\Z/2)^k$. 
Since it is non-degenerate, the correspondence 
\begin{equation} \label{hom}
(\Z/2)^k\to \Hom((\Z/2)^k,\Z/2) \quad\text{given by $u\to (u,\ )$}
\end{equation}
is an isomorphism, where $\Hom((\Z/2)^k,\Z/2)$ denotes the group of homomorphisms from 
$(\Z/2)^k$ to $\Z/2$.

\begin{lemm} \label{lemm3}
If $u\in (\Z/2)^k$ is non-zero, then the kernel of $(u,\ )$ is a codimension 1 subspace 
of $(\Z/2)^k$.  On the other hand, 
any codimension 1 subspace $V$ of $(\Z/2)^k$ is obtained as the kernel of $(u,\ )$ 
for some non-zero $u\in (\Z/2)^k$ and $u$ is uniquely determined by $V$.
\end{lemm}

\begin{proof} 
The former statement in the lemma follows from the fact that the bilinear form $(\ ,\ )$ is non-degenerate. 
Let $V$ be a codimension 1 subspace of $(\Z/2)^k$.  Then the quotient vector space of $(\Z/2)^k$ by 
$V$ is one-dimensional, so it is isomorphic to $\Z/2$ and hence defines an element of $\Hom((\Z/2)^k,\Z/2)$ 
whose kernel is $V$.  This together with \eqref{hom} implies the latter statement in the lemma. 
\end{proof}

\begin{lemm} \label{prop3}
Suppose $k\geqq 2$.  Then 
$s_\R(m,p)\geqq k$ if and only if there is a set of non-negative integers $\{ a_v\mid v\in (\Z/2)^k\backslash\{0\}\}$ 
with $\sum a_v=m$,  
which satisfy the following $(2^k-1)$ inequalities 
\begin{equation*} \label{eq2'}
\sum_{(u,v)=0} a_v\leqq p-1\quad \text{for each $u\in (\Z/2)^k\backslash\{0\}$}.
\end{equation*}
\end{lemm}

\begin{proof}
Any codimension $1$ subspace of $(\Z/2)^k$ is the kernel of 
a homomorphism $(u,\ )\colon (\Z/2)^k\to \Z/2$ for some non-zero $u\in (\Z/2)^k$ by Lemma~\ref{lemm3}. 
Therefore any $p$ column vectors in a $k\times m$ matrix with $a_v$ numbers of column vector $v$ 
for each $v$ span $(\Z/2)^k$ if and only if the $a_v$'s satisfy the inequalities in the lemma.  This proves the lemma. 
\end{proof}

The lemma above shows that our problem is a problem of \emph{integer} linear programming.  
If we consider the problem over real numbers, then it is easy to find the solution of the 
problem as shown by the following lemma. 

\begin{lemm} \label{lemm:max}
Suppose that $k\geqq 2$  and let $b$ be a real number.  
If we allow $a_v$'s to be real numbers and $a_v$'s satisfy the following $(2^k-1)$ inequalities 
\begin{equation} \label{eq2}
\sum_{(u,v)=0} a_v\leqq  b\quad \text{for each $u\in (\Z/2)^k\backslash\{0\}$}, 
\end{equation} 
then the linear function $\sum a_v$ on $\R^{2^k-1}$ takes the maximum value $$(2^k-1)b/(2^{k-1}-1)$$ at 
a unique point $x=(a_v)\in\R^{2^k-1}$ with $a_v=b/(2^{k-1}-1)$ for every $v$. 
\end{lemm}

\begin{proof}
Each $a_v$ appears in exactly $(2^{k-1}-1)$ times in the inequalities \eqref{eq2} because 
there are exactly $(2^{k-1}-1)$ numbers of $u\in (\Z/2)^k\backslash\{0\}$ such that $(u,v)=0$.  
Therefore, taking sum of the $(2^k-1)$ inequalities \eqref{eq2} over $u\in (\Z/2)^k\backslash\{0\}$, we obtain 
\[
(2^{k-1}-1)\sum a_v\leqq (2^k-1)b
\]
and the equality is attained at the point $x$ in the lemma; so the maximum value of $\sum a_v$ 
satisfying \eqref{eq2} is $(2^k-1)b/(2^{k-1}-1)$.  

We shall observe that the maximum value $(2^k-1)b/(2^{k-1}-1)$ is attained only at the point $x$. 
Suppose that $\sum a_v$ takes the maximum value on $a_v$'s satisfying \eqref{eq2}.  
Then the argument above shows that all the inequalities in \eqref{eq2} 
must be equalities, i.e. 
\begin{equation} \label{eqb}
\sum_{(u,v)=0}a_v=b \quad \text{for each $u\in (\Z/2)^k\backslash\{0\}$}. 
\end{equation}
We choose one $v$ arbitrarily and take sum of \eqref{eqb} over all non-zero $u$'s with $(u,v)=0$.  
The number of such $u$ is $2^{k-1}-1$, so $a_v$ appears $2^{k-1}-1$ times in the sum.  
But $a_{v'}$ with $v'\not=v$ appears $2^{k-2}-1$ times in the sum because the number of non-zero $u$ with 
$(u,v)=(u,v')=0$ is $2^{k-2}-1$.  Therefore we obtain
\begin{equation} \label{eqb2}
(2^{k-1}-1)a_v+(2^{k-2}-1)\sum_{v'\not=v}a_{v'}=(2^{k-1}-1)b.
\end{equation}
Here 
\begin{equation} \label{eqb3}
\sum_{v'\not=v}a_{v'}=(2^k-1)b/(2^{k-1}-1)-a_v
\end{equation}
since $\sum_v a_v$ is assumed to take the maximum value $(2^k-1)b/(2^{k-1}-1)$.  
Plugging \eqref{eqb3} in \eqref{eqb2}, we obtain 
\[
2^{k-2}a_v+(2^{k-2}-1)\frac{(2^k-1)b}{(2^{k-1}-1)}=(2^{k-1}-1)b
\]
and a simple computation shows $a_v=b/(2^{k-1}-1)$. 
\end{proof}

Lemma~\ref{lemm:max} tells us that the point $x$ is a unique vertex of the 
polyhedron $P(b)$ defined by the inequalities \eqref{eq2} and $(2^k-1)$ hyperplanes 
$\sum_{(u,v)=0} a_v=b$ in $\R^{2^k-1}$ $(u\in (\Z/2)^k\backslash\{0\})$ are 
in general position.  
Motivated by Lemma~\ref{prop3} we make the following definition. 

\begin{defi}
For a positive integer $k\geqq 2$ and a non-negative integer $b$, 
we define $m_k(b)$ to be the maximum integer which the linear function $\sum a_v$ 
takes on lattice points satisfying \eqref{eq2} and $a_v\geqq 0$ for every $v$.  
\end{defi}

One easily sees that $m_k(0)=0$ and $m_k(b)\geqq b$ for any $b$. 
The importance of finding values of $m_k(b)$ lies in the following lemma. 

\begin{lemm} \label{lemm:sR}
$s_\R(m,p)=k$ for $k\geqq 2$ if and only if $m_{k+1}(p-1)<m\leqq m_k(p-1)$.  
\end{lemm}

\begin{proof}
Since $s_\R(m,p)$ decreases as $m$ increases by Proposition~\ref{prop2} (3), the lemma follows from 
Lemma~\ref{prop3}.
\end{proof} 

\begin{rema} \label{rema1}
Since $s_\R(m,p)\leqq p$ by Proposition~\ref{prop2} (1), the equality $s_\R(m,p)=k$ makes sense only when $k\leqq p$.  
In other words, $m_k(b)$ has the matrix interpretation discussed for $s_\R(m,p)$ in Section~\ref{sect:3} only when $k\leqq b+1$.  
\end{rema}

The following is essentially a restatement of Theorem~\ref{mm+1}.  

\begin{theo} \label{mod2}
$m_k(b)\equiv b \pmod 2$.
\end{theo}

\begin{proof}
It is not difficult to see that $m_k(b)=b$ when $b\leqq k-2$ (see Theorem~\ref{mkk}), so the the theorem 
holds in this case.  Suppose $b\geqq k-1$ and set $b=p-1$.  Then $s_\R(m_k(p-1),p)=k$ by Lemma~\ref{lemm:sR}.  
If $m_k(p-1)-p$ is even, then $s_\R(m_k(p-1)+1,p)=k$ by Theorem~\ref{mm+1}.  But this contradicts 
the maximality of $m_k(p-1)$.  Therefore $m_k(p-1)-p$ is odd, i.e., $m_k(b)-b$ is even. 
\end{proof}

The following corollary follows from Lemma~\ref{lemm:max} 
and the last statement in the corollary also follows from Theorem~\ref{theo3}. 

\begin{coro} \label{coro:divisi}
For any non-negative integer $b$ we have 
\begin{equation} \label{eq:upper}
m_k(b)\leqq \Big[ \frac{(2^k-1)b}{2^{k-1}-1}\Big]=2b+\Big[\frac{b}{2^{k-1}-1}\Big]
\end{equation}
and the equality is attained when $b$ is divisible by $2^{k-1}-1$, i.e. 
\begin{equation} \label{eq:divisi}
\text{$m_k((2^{k-1}-1)Q)=(2^k-1)Q$}
\end{equation}
for any non-negative integer $Q$. In particular 
\begin{equation} \label{eq:m2b}
\text{$m_2(b)=3b$ for any $b$.}
\end{equation} 
\end{coro} 

One can find some values of $s_\R(m,p)$ using \eqref{eq:divisi}.  

\begin{exam}
Take $p=(2^{k-1}-1)(2^k-1)q+1$ where $q$ is any positive integer.  Then
\[
\text{$m_{k}(p-1)=(2^{k}-1)^2q$,\quad $m_{k+1}(p-1)=(2^{k+1}-1)(2^{k-1}-1)q$}
\]  
by \eqref{eq:divisi}.  
Therefore it follows from Lemma~\ref{lemm:sR} that 
$s_R(m,p)=k$ for $m$ with $(2^{k+1}-1)(2^{k-1}-1)q<m\leqq (2^{k}-1)^2q$.  
\end{exam}

\section{Some more properties of $m_k(b)$} \label{sect:5}

In this section, we study some more properties of $m_k(b)$.  

\begin{lemm} \label{coro3}   For any non-negative integers $b,b'$ we have 
\begin{equation} \label{eq:additive}
m_k(b)+m_k(b')\leqq m_k(b+b').
\end{equation}
In particular, 
\begin{enumerate}
\item $m_k(b)+b'\leqq m_k(b+b')$, 
\item $m_k(b)+(2^k-1)Q\leqq m_k(b+(2^{k-1}-1)Q)$ for any non-negative integer $Q$.
\end{enumerate}
\end{lemm}

\begin{proof}
Let $\{a_v\}$ (resp. $\{a'_v\}$) be a set of non-negative integers which satisfy \eqref{eq2} and $\sum a_v=m_k(b)$ 
(resp. \eqref{eq2} with $b$ replaced by $b'$ and $\sum a'_v=m_k(b')$).  Then $\{a_v+a'_v\}$ is a set of non-negative integers which 
satisfy \eqref{eq2} with $b$ replaced by $b+b'$ and $\sum (a_v+a'_v)=m_k(b)+m_k(b')$.  Therefore \eqref{eq:additive} follows.   

The inequality (1) follows from \eqref{eq:additive} and the fact that $m_k(b')\geqq b'$.  
The inequality (2) follows by taking $b'=(2^{k-1}-1)Q$ in \eqref{eq:additive} and using \eqref{eq:divisi}. 
\end{proof}

We will see in later sections that the equality in Lemma~\ref{coro3} (1)  holds for special values of $b$ and $b'$ but 
does not hold in general.  However, \eqref{eq:divisi} and results obtained in later sections imply that 
the equality in Lemma~\ref{coro3} (2) would hold for arbitrary values of $b$ and $Q$.
We shall formulate it as the following conjecture.

\begin{conj}
$m_k((2^{k-1}-1)Q+R)=(2^k-1)Q+m_k(R)$ for any non-negative integers $Q$ and $R$, where we may assume 
$0\leqq R\leqq 2^{k-1}-2$ without loss of generality. 
\end{conj}

The following lemma enables us to find an upper bound for $m_k(b)$ by induction on $k$ and we will see 
that the former inequality in \eqref{eq:mkbk-1} is not always but often an equality.  

\begin{lemm} \label{lemm:kk-1}
If $b$ is not divisible by $2^{k-1}-1$ and $Q=[b/(2^{k-1}-1)]$, then  
\[
m_k(b)\leqq m_{k-1}(b-q-1)+q+1
\]
for any integer $0\leqq q\leqq Q$ and $m_{k-1}(b-q-1)+q+1$ increases as $q$ decreases; so in particular 
\begin{equation} \label{eq:mkbk-1}
m_k(b)\leqq m_{k-1}(b-Q-1)+Q+1\leqq m_{k-1}(b-1)+1.
\end{equation}
\end{lemm}

\begin{proof}
Let $\{a_v\}$ be a set of non-negative integers which satisfy \eqref{eq2} and $\sum a_v=m_k(b)$.  
Then 
\begin{equation} \label{uv}
\text{$\sum_{(u,v)=0}a_v=b$ for some $u\in (\Z/2)^k\backslash\{0\}$}
\end{equation}
because otherwise we can add $1$ to 
some $a_v$ so that the resulting set of non-negative integers still satisfy \eqref{eq2} but their sum is $m_k(b)+1$, which 
contradicts the definition of $m_k(b)$.  Therefore $a_{v}\geqq Q+1$ for some $a_v$ in \eqref{uv} because if $a_v\leqq Q$ for any $v$, 
then $\sum_{(u,v)=0} a_v\leqq (2^{k-1}-1)Q$ and $(2^{k-1}-1)Q$ 
is strictly smaller than $b$ since $b$ is not divisible by $2^{k-1}-1$ by assumption.  

Through a linear transformation of $(\Z/2)^k$, 
we may assume that the $v$ with $a_v\geqq Q+1$ is $\e_k=(0,\dots,0,1)^T$, so
\begin{equation} \label{aek}
a_{\e_k}\geqq Q+1.
\end{equation}
The kernel $\e_k^\perp$ of the homomorphism $(\e_k,\ ): (\Z/2)^k\to \Z/2$ 
can naturally be identified with $(\Z/2)^{k-1}$.  For $u\in \e_k^\perp$, \eqref{eq2} reduces to 
\begin{equation} \label{ak}
a_{\e_k}+\sum_{(u,v)=0, v\not=\e_k}a_v\leqq b.
\end{equation}
Let $\pi\colon (\Z/2)^k\to (\Z/2)^{k-1}$ be the natural projection.  For $u\in \e_k^\perp$, we have $(u,v)=0$ if and only if 
$(\pi(u),\pi(v))=0$.  Therefore \eqref{ak} reduces to 
\[
\sum_{(\pi(u),\bar v)=0}a_{\bar v}\leqq b-a_{\e_k}
\]
where $\bar v$ runs over all non-zero elements of $(\Z/2)^{k-1}$ and $a_{\bar v}=\sum_{\pi(v)=\bar v} a_v$.    
It follows that $\sum a_{\bar v}\leqq m_{k-1}(b-a_{\e_k})$ and hence 
\begin{equation} \label{eq:mkb}
m_k(b)=\sum a_v=a_{\e_k}+\sum a_{\bar v}\leqq a_{\e_k}+m_{k-1}(b-a_{\e_k}).
\end{equation}
Here $q+m_{k-1}(b-q)$ increases as $q$ decreases because it follows from Lemma~\ref{coro3} that 
\[
q+m_{k-1}(b-q)\leqq q-1+m_{k-1}(b-q+1). 
\]
Therefore, the inequalities in the lemma follow from \eqref{eq:mkb} and \eqref{aek}. 
\end{proof}

\begin{coro} \label{coro:mdec}
$m_k(b)\leqq m_{k-1}(b)$ for any $b$ and $k\geqq 3$. 
\end{coro}

\begin{proof}
Since $m_{k-1}(b-q-1)+q+1\leqq m_{k-1}(b)$ by Lemma~\ref{coro3} (1), the corollary follows from Lemma~\ref{lemm:kk-1}. 
\end{proof}

We shall give another application of Lemma~\ref{lemm:kk-1}.  
Our conjecture stated in this section can be thought of as a periodicity of $m_k(b)$ for a fixed $k$. 
The following proposition implies another periodicity of $m_k(b)$, where $k$ varies. 
It in particular says that once we know values 
of $m_k(b)$ for all $b$, we can find values of $m_{k+1}(b)$ for \lq\lq half" of all $b$.    

\begin{prop} \label{period}
Suppose that  $$m_k((2^{k-1}-1)Q+R)=(2^k-1)Q+m_k(R)$$ for some $k, R$ and any $Q$ where $0\leqq R\leqq 2^{k-1}-2$.  
Then 
\begin{equation} \label{f}
m_{k+1}((2^{k}-1)Q+2^{k-1}+R)=(2^{k+1}-1)Q+2^k+m_k(R),
\end{equation} 
more generally, 
\begin{equation} \label{l}
m_{k+\ell}((2^{k+\ell-1}-1)Q+2^{k+\ell-1}-2^{k-1}+R)=(2^{k+\ell}-1)Q+2^{k+\ell}-2^k+m_k(R)
\end{equation}
for any non-negative integer $\ell$. 
\end{prop}  

\begin{proof}
The latter identity \eqref{l} easily follows if we use the former statement repeatedly,  so we prove only \eqref{f}.  
When $R=0$,  \eqref{f} follows from \eqref{eq:divisi}; so we may assume $R\not=0$. 
It follows from Lemma~\ref{lemm:kk-1} and the assumption in the lemma that
\begin{equation} 
\begin{split}
&m_{k+1}((2^k-1)Q+2^{k-1}+R)\\
\leqq& m_k((2^k-1)Q+2^{k-1}+R-Q-1)+Q+1\\
=&m_k((2^{k-1}-1)(2Q+1)+R)+Q+1\\
=&(2^k-1)(2Q+1)+m_k(R)+Q+1\\
=&(2^{k+1}-1)Q+2^k+m_k(R).
\end{split}
\end{equation}

We shall prove the opposite inequality.  
Let $\{a_v\}$ be a set of non-negative integers which satisfy \eqref{eq2} with $b$ replaced by $R$ and 
\begin{equation} \label{sumav}
\text{$\sum a_v=m_k(R)$.}
\end{equation}  
We regard $(\Z/2)^k$ as a subspace of $(\Z/2)^{k+1}$ in a natural way and define 
$a_v'$ for $v\in (\Z/2)^{k+1}$ by 
\begin{equation} \label{av'}
a'_v:=\begin{cases} Q+a_v \quad&\text{for $v\in (\Z/2)^k\backslash\{0\}$,}\\
Q+1 \quad&\text{for $v\notin (\Z/2)^k$.}
\end{cases}
\end{equation}
We shall check that 
the set $\{a'_v\}$ of non-negative integers satisfies \eqref{eq2} with $b$ replaced by 
\begin{equation} \label{b'}
\text{$b':=(2^k-1)Q+2^{k-1}+R$.}
\end{equation}   
Let $u\in (\Z/2)^{k+1}\backslash\{0\}$ and denote by $u^\perp$ the kernel of the homomorphism 
$(u,\ )\colon (\Z/2)^{k+1}\to \Z/2$, which is a codimension 1 subspace of $(\Z/2)^{k+1}$.  
We distinguish two cases. 

{\bf Case 1.} The case where $u^\perp=(\Z/2)^k$. 
It follows from \eqref{sumav} and \eqref{av'} that 
\begin{equation} \label{sumav'}
\begin{split}
\sum_{(u,v)=0}a'_v&=\sum (Q+a_v)\\
&=(2^{k}-1)Q+\sum a_v\\
&=(2^k-1)Q+m_k(R).
\end{split}
\end{equation}
Here $m_k(R)\leqq 2R$ by \eqref{eq:upper} and since $R\leqq 2^{k-1}-2$, we obtain
\begin{equation*} \label{eqmkR}
\text{$m_k(R)\leqq 2^{k-1}+R$.}
\end{equation*} 
This together with \eqref{b'} and \eqref{sumav'} shows that $\sum_{(u,v)=0}a_v'\leqq b'$.  

{\bf Case 2.}  The case where $u^\perp\not=(\Z/2)^k$.  Since both $u^\perp$ and $(\Z/2)^k$ are codimension 1 
subspaces of $(\Z/2)^{k+1}$ and they are different, the intersection $u^\perp\cap (\Z/2)^k$ is a codimension 1 
subspace of $(\Z/2)^k$ and hence the number of elements in $u^\perp\backslash (\Z/2)^k$ is $2^{k-1}$.  
Therefore, it follows from \eqref{av'} and \eqref{b'} that 
\[
\begin{split}
\sum_{(u,v)=0}a_v'&=\sum_{v\in u^\perp\cap(\Z/2)^k}a_v'+\sum_{v\in u^\perp\backslash (\Z/2)^k}a_v'\\
&=\sum_{v\in u^\perp\cap(\Z/2)^k}(Q+a_v)+\sum_{v\in u^\perp\backslash (\Z/2)^k}(Q+1)\\
&= (2^k-1)Q+\sum_{v\in u^\perp\cap(\Z/2)^k}a_v +2^{k-1}\\
&\leqq (2^k-1)Q+R+2^{k-1}=b'
\end{split}
\]
where the inequality above follows from the fact that the set $\{a_v\}$ satisfies \eqref{eq2} with $b$ replaced by $R$. 

The above two cases prove that the set $\{a'_v\}$ satisfies \eqref{eq2} with $b$ replaced by $b'$.  Finally it follows from 
\eqref{sumav} and  \eqref{av'} that  
\[
\begin{split}
\sum_{v\in(\Z/2)^{k+1}\backslash\{0\}} a_v'&=\sum_{v\in(\Z/2)^k\backslash\{0\}} (Q+a_v) +\sum_{v\notin (\Z/2)^k}(Q+1)\\
&=(2^{k+1}-1)Q+\sum_{v\in(\Z/2)^k\backslash\{0\}} a_v+2^k\\
&=(2^{k+1}-1)Q+m_k(R)+2^k.
\end{split}
\]
This implies the following desired opposite inequality 
$$m_{k+1}((2^k-1)Q+2^{k-1}+R)\geqq (2^{k+1}-1)Q+2^k+m_k(R)$$
and completes the proof of \eqref{f}. 
\end{proof}

\section{$m_k(b)$ for $b\leqq k+1$} \label{sect:6}

In this section we will find the values of $m_k(b)$ for $b\leqq k+1$.  
We treat the case where $b\leqq k-1$ first. 

\begin{theo} \label{mkk}
For any $k\geqq 2$, we have 
\[
m_k(b)=\begin{cases} b \quad&\text{if $b\leqq k-2$,}\\
b+2 \quad&\text{if $b=k-1$}. \end{cases}
\]
\end{theo}

\begin{proof}
(1) The case where $b\leqq k-2$.  Let $a_v$'s be non-negative integers which 
satisfy \eqref{eq2}.  Suppose that there are more than $b$ positive integers  $a_v$'s and 
choose $b+1$ out of them.  Since $b+1\leqq k-1$, $v$'s for the chosen $b+1$ positive $a_v$'s are contained 
in some codimension 1 subspace of $(\Z/2)^k$; so the sum of those $b+1$ positive $a_v$'s must be 
less than or equal to $b$ by \eqref{eq2}, which is a contradiction.  
Therefore there are at most $b$ positive $a_v$'s.  Since $b\leqq k-2$, $v$'s for the positive $a_v$'s 
are contained in some codimension 1 subspace of $(\Z/2)^k$; so $\sum a_v\leqq b$ by \eqref{eq2} 
and this proves $m_k(b)\leqq b$.  On the other hand, it is clear that $m_k(b)\geqq b$, so 
$m_k(b)=b$ when $b\leqq k-2$. 

(2) The case where $b=k-1$.  In this case we can use the matrix interpretation of $m_k(b)$, 
see the Remark following Lemma~\ref{lemm:sR}. 
The following argument is essentially same as Lemma~\ref{lemm2}.  
Let $A$ be a $k\times m$ matrix where any $k$ column vectors span $(\Z/2)^k$.  
We may assume that the first $k$ column vectors are the standard basis, so 
$A=(\e_1,\dots,\e_k,\a_{k+1},\dots,\a_m)$.  Since any $k-1$ vectors from $\e_1,\dots,\e_k$ 
together with $\a_{j}$ span $(\Z/2)^k$, $\a_{j}$ must be $\sum_{i=1}^k \e_i$.  
Therefore $m$ must be less than or equal to $k+1$ and this shows $m_k(k-1)\leqq k+1$. 
On the other hand, since any $k$ column vectors in $(\e_1,\dots,\e_k,\sum\e_i)$ span 
$(\Z/2)^k$, $m_k(k-1)\geqq k+1$.  This proves $m_k(k-1)=k+1$. 
\end{proof}

\begin{theo} \label{mkk0}
If $b=k$, then 
\[
m_k(b)=\begin{cases} b+4 \quad&\text{if $k=2,3,4$},\\
b+2\quad&\text{if $k\geqq 5$}.
\end{cases}
\]
\end{theo}

\begin{proof}
Since $m_2(2)=6$ by \eqref{eq:m2b} and $m_3(3)=7$ by \eqref{eq:divisi}, the theorem is proven when $k=2,3$. 
One can easily check that any 5 columns in this matrix 
\[
{\tiny \begin{pmatrix} 1&0&0&0&0&1&1&1\\
0&1&0&0&1&0&1&1\\
0&0&1&0&1&1&0&1\\
0&0&0&1&1&1&1&0\end{pmatrix}}
\]
span $(\Z/2)^4$, so $m_4(4)\geqq 8$. 
On the other hand, using Lemma~\ref{lemm:kk-1}, we obtain 
\[
m_4(4)\leqq m_{3}(3)+1=8.  
\]
Thus $m_4(4)=8$ and the theorem is proven when $k=4$. 

Since $m_k(k-1)=k+1$ by Theorem~\ref{mkk}, it follows from  Lemma~\ref{coro3} (1) that 
\[
m_k(k)\geqq m_k(k-1)+1=k+2.
\]
In the sequel it suffices to prove that if $m_k(k)\geqq k+3$, then $k\leqq 4$.  

Suppose $m_k(k)\geqq k+3$.  Then there is a $k\times (k+3)$ matrix $A$ with entries in $\Z/2$ 
such that any $k+1$ column vectors in $A$ span $(\Z/2)^k$.  We may assume that 
$A=(\e_1,\dots,\e_k,\a_1,\a_2,\a_3)$ as before. 
Denote by $\a^i$ the $i$-th row vector in the submatrix $(\a_1,\a_2,\a_3)$.    
Since any $k+1$ column vectors in $A$ span $(\Z/2)^k$, we see that 
\[
\begin{pmatrix} \a^i\\
\a^j \end{pmatrix}=\begin{pmatrix} 1&0&1\\
0&1&1\end{pmatrix},\quad \begin{pmatrix} 1&1&1\\
0&1&1\end{pmatrix}\quad\text{or}\quad \begin{pmatrix} 0&1&1\\
1&1&1\end{pmatrix}
\]
up to permutations of column vectors at the right hand side.   This must occur for any $1\leqq i<j\leqq k$ 
but one can easily see that this is impossible when $k\geqq 5$.   
\end{proof}

\begin{theo} \label{mkk1}
If $b=k+1$, then 
\[
m_k(b)=\begin{cases} b+6 \quad&\text{if $k=2$},\\
b+4 \quad&\text{if $3\leqq k\leqq 11$},\\
b+2\quad&\text{if $k\geqq 12$}.
\end{cases}
\]
\end{theo}

\begin{proof}
Since $m_2(3)=9$ by \eqref{eq:m2b}, the theorem is proven when $k=2$.   

Using Lemma~\ref{lemm:kk-1} repeatedly, we have
\begin{equation} \label{m11}
m_{11}(12)\leqq m_{10}(11)+1\leqq m_9(10)+2\leqq \dots\leqq m_3(4)+8\leqq m_2(2)+10=16
\end{equation}
where we used \eqref{eq:m2b} at the last identity.  
On the other hand, it follows from Theorem~\ref{theo2} that 
\[
s_\R(16,13)=s_\R(15,13)=[15-\log_2(15+1)]=11
\]
and hence  $m_{11}(12)\geqq 16$ by Lemma~\ref{lemm:sR}.  Therefore 
$m_{11}(12)=16$ and all the inequalities in \eqref{m11} must be equalities, proving the 
second case in the theorem. 

Similarly, it follows from Theorem~\ref{theo2}  that 
$$s_\R(16,14)=[16-\log_2(16+1)]=11$$
and hence 
$m_{12}(13)\leqq 15$ 
by Lemma~\ref{lemm:sR}.  On the other hand, it follows from Theorem~\ref{mkk0} and Corollary~\ref{coro:mdec} that 
\[
15=m_{13}(13)\leqq m_{12}(13).
\]
Therefore $m_{12}(13)=15$. 

Suppose $k\geqq 12$.  Then using Lemma~\ref{lemm:kk-1} repeatedly, we have 
\[
m_k(k+1)\leqq m_{k-1}(k)+1\leqq \dots\leqq m_{12}(13)+k-12=k+3
\]
where we used the fact $m_{12}(13)=15$ just shown above. 
On the other hand, it follows from Lemma~\ref{coro3} (1) and 
Theorem~\ref{mkk0} that 
\[
m_k(k+1)\geqq m_k(k)+1=k+3.
\]
Therefore $m_k(k+1)=k+3$ when $k\geqq 12$, 
proving the last case in the theorem. 
\end{proof}

\section{Further computations of $m_k(b)$} \label{sect:56}

In this section we will make some more computations of $m_k(b)$ by combining the results in the previous sections.  
All of the results provide supporting evidence to the Conjecture stated in Section~\ref{sect:5}. 

\begin{prop} \label{theo:Rk-1}
If $R\leqq k-1$, then 
$$m_k((2^{k-1}-1)Q+R)=(2^k-1)Q+m_k(R)$$
where 
\[
m_k(R)=\begin{cases} R \quad&\text{if $R\leqq k-2$,}\\
R+2 \quad&\text{if $R=k-1$}. \end{cases}
\]
by Theorem~\ref{mkk}. 
\end{prop}

\begin{proof}
When $R=0$, the proposition follows from \eqref{eq:divisi} since $m_k(0)=0$.  
So we may assume $1\leqq R\leqq k-1$.  We prove the proposition by induction on $k$. 
Since $m_2(b)=3b$ by \eqref{eq:m2b}, the proposition holds when $k=2$.  
Suppose the proposition holds for $k=\ell-1$.  
It follows from \eqref{eq:divisi}, Lemmas~\ref{coro3}, 
~\ref{lemm:kk-1} and the induction assumption that 
\begin{equation} \label{2k-1}
\begin{split}
(2^{\ell}-1)Q+m_{\ell}(R)&=m_{\ell}((2^{\ell-1}-1)Q)+m_{\ell}(R)\\
&\leqq m_{\ell}((2^{\ell-1}-1)Q+R)\\
&\leqq m_{\ell-1}((2^{\ell-1}-1)Q+R-Q-1)+Q+1\\
&=m_{\ell-1}((2^{\ell-2}-1)2Q+R-1)+Q+1\\
&=(2^{\ell-1}-1)2Q+m_{\ell-1}(R-1)+Q+1\\
&=(2^{\ell}-1)Q+m_{\ell-1}(R-1)+1.
\end{split}
\end{equation}
Here since $R\leqq \ell-1$, we have 
$m_\ell(R)=m_{\ell-1}(R-1)+1$ 
by Theorem~\ref{mkk}.  Therefore the first and last terms in \eqref{2k-1} are same, 
so the first inequality in \eqref{2k-1} must be an equality, which proves the proposition
when $k=\ell$, completing the induction step.  
\end{proof}

The following corollary follows from Proposition~\ref{theo:Rk-1} by taking $k=3$. 

\begin{coro} \label{m3b}
\[
m_3(3Q+R)=\begin{cases} 7Q \quad&\text{if $R=0$},\\
7Q+1\quad&\text{if $R=1$},\\
7Q+4 \quad&\text{if $R=2$}.
\end{cases}
\]
\end{coro}

Combining Proposition~\ref{theo:Rk-1} with Proposition~\ref{period}, 
one can improve Proposition~\ref{theo:Rk-1} as follows.   

\begin{theo} \label{theo:R=1}
Let $0\leqq \ell\leqq k-2$.  If $0\leqq r\leqq k-\ell-1$, then 
\[
m_k((2^{k-1}-1)Q+2^{k-1}-2^{k-1-\ell}+r)=(2^k-1)Q+2^k-2^{k-\ell}+m_{k-\ell}(r)
\]
where 
\[
m_{k-\ell}(r)=\begin{cases} r \quad&\text{if $r\leqq k-\ell-2$,}\\
r+2 \quad&\text{if $r=k-\ell-1$}. \end{cases}
\]
by Theorem~\ref{mkk}. 
\end{theo}

\begin{proof}
By Proposition~\ref{theo:Rk-1}, we have 
\begin{equation} \label{condr}
m_k((2^{k-1}-1)Q+r)=(2^k-1)=(2^k-1)Q+m_k(r) \quad \text{for $0\leqq r\leqq k-1$}.
\end{equation}
Therefore, it follows from \eqref{l} in Proposition~\ref{period} that 
\begin{equation} \label{mkell}
m_{k+\ell}((2^{k+\ell-1}-1)Q+2^{k+\ell-1}-2^{k-1}+r)=(2^{k+\ell}-1)Q+2^{k+\ell}-2^k+m_k(r)
\end{equation}
for any non-negative integer $\ell$.  Rewriting $k+\ell$ as $k$, the identity \eqref{mkell} 
turns into the identity in the theorem and the condition $0\leqq r\leqq k-1$ in \eqref{condr} turns into 
the condition $0\leqq r\leqq k-\ell-1$ in the theorem. 
\end{proof}

\begin{prop} \label{prop:Rk+1}
If $R=k+1$ and $4\leqq k\leqq 11$, then 
\[
m_k((2^{k-1}-1)Q+R)=(2^k-1)Q+m_k(R)
\]
where $m_k(R)=R+4$ by Theorem~\ref{mkk1}. 
\end{prop}

\begin{proof}
First we prove the proposition when $k=4$.  In this case $R=5$. 
It follows from Lemma~\ref{lemm:kk-1} and Corollary~\ref{m3b} that 
\[
\begin{split}
m_4((2^3-1)Q+5)&\leqq m_3(7Q+5-Q-1)+Q+1\\
&=7(2Q+1)+1+Q+1=15Q+9
\end{split}
\]
while it follows from \eqref{eq:additive}, \eqref{eq:divisi} and Theorem~\ref{mkk1} 
\[
\begin{split}
m_4((2^3-1)Q+5)&\geqq m_4((2^3-1)Q)+m_4(5)\\
&=(2^4-1)Q+9=15Q+9.
\end{split}
\]
This proves the proposition when $k=4$.  

Suppose that the proposition holds for $k-1$ with $4\leqq k-1\leqq 10$.  
Then it follows from Lemma~\ref{lemm:kk-1} and the induction assumption that 
\[
\begin{split}
m_k((2^{k-1}-1)Q+R)&\leqq m_{k-1}((2^{k-1}-1)Q+R-Q-1)+Q+1\\
&=m_{k-1}((2^{k-2}-1)2Q+R-1)+Q+1\\
&=(2^{k-1}-1)2Q+(R-1)+4+Q+1\\
&=(2^k-1)Q+R+4
\end{split}
\]
while it follows from \eqref{eq:additive}, \eqref{eq:divisi} and Theorem~\ref{mkk1} 
\[
\begin{split}
m_k((2^{k-1}-1)Q+R)&\geqq m_k((2^{k-1}-1)Q)+m_k(R)\\
&=(2^k-1)Q+R+4.
\end{split}
\]
These show that $m_k((2^{k-1}-1)Q+R)=(2^k-1)Q+R+4$, completing the induction step.    
\end{proof}

Similarly to Theorem~\ref{theo:R=1}, Proposition~\ref{prop:Rk+1} can be improved as follows by 
combining it with Proposition~\ref{period}.  The proof is same as that of 
Theorem~\ref{theo:R=1}, so we omit it.

\begin{theo} \label{theo:Rk+1}
Let $0\leqq \ell\leqq k-2$.  If $4\leqq k-\ell\leqq 11$, then  
\[
m_k((2^{k-1}-1)Q+2^{k-1}-2^{k-\ell-1}+k-\ell+1)=(2^k-1)Q+2^k-2^{k-\ell}+k-\ell+5.
\]
\end{theo}

\begin{exam}  \label{ex:k=34}
Below is a table of values of $m_k((2^{k-1}-1)Q+R)$ 
for $k=2,3,4,5,6$.   
{\tiny
\begin{table}[htb]
\begin{tabular}{|c||c|c|c|c|c|} \hline
$R\backslash k$& 2 &  3 & 4 & 5 & 6  \\
\hline\hline 0 & 3Q & 7Q    &  15Q  &   31Q           &    63Q           \\
\hline  1        &  & 7Q+1 & 15Q+1 &  31Q+1           &    63Q+1       \\
\hline  2        &  & 7Q+4 & 15Q+2 &  31Q+2           &    63Q+2       \\
\hline  3        &  &          & 15Q+5 &  31Q+3           &    63Q+3       \\
\hline 4         &  &          & 15Q+8 &  31Q+6           &    63Q+4       \\
\hline 5         &  &          & 15Q+9 &  31Q+7 or 9    &    63Q+7       \\
\hline 6         &  &         & 15Q+12 & 31Q+10         &    63Q+8 or 10 \\
\hline 7         &  &         &             &  31Q+11 or 13&    63Q+11      \\
\hline 8         &  &         &             &  31Q+16         &    63Q+12 or 14 \\
\hline 9         &  &         &             &  31Q+17         &    63Q+13, 15 or 17 \\
\hline 10       &  &         &             &  31Q+18         &    63Q+14, 16 or 18 \\
\hline 11       &  &         &             &  31Q+21         &    63Q+15, 17 or 19 \\
\hline 12       &  &         &             &  31Q+24         &    63Q+20 or 22 \\
\hline 13       &  &         &             &  31Q+25         &    63Q+21, 23 or 25\\
\hline 14       &  &         &             &  31Q+28         &    63Q+24 or 26  \\
\hline 15       &  &         &             &                      &     63Q+27 or 29\\
\hline 16       &  &         &             &                      &     63Q+32     \\
\hline 17       &  &         &             &                      &     63Q+33     \\
\hline 18       &  &         &             &                      &     63Q+34     \\
\hline 19       &  &         &             &                      &     63Q+35     \\
\hline 20       &  &         &             &                      &     63Q+38     \\
\hline 21       &  &         &             &                      &     63Q+39 or 41\\
\hline 22       &  &         &             &                      &     63Q+42     \\
\hline 23       &  &         &             &                      &     63Q+43 or 45  \\
\hline 24       &  &         &             &                      &     63Q+48     \\
\hline 25       &  &         &             &                      &     63Q+49     \\
\hline 26       &  &         &             &                      &     63Q+50     \\
\hline 27       &  &         &             &                      &     63Q+53     \\
\hline 28       &  &         &             &                      &     63Q+56     \\
\hline 29       &  &         &             &                      &     63Q+57     \\
\hline 30       &  &         &             &                      &     63Q+60     \\
\hline
\end{tabular} \vspace{1.25ex}
\caption{$m_k((2^{k-1}-1)Q+R)$ for $k=3,4,5,6$}
\label{table1}
\end{table}
}

The values above for $k=2, 3, 4$ can be obtained from Theorem~\ref{theo:R=1} although they are obtained from   
\eqref{eq:m2b} when $k=2$ and from Corollary~\ref{m3b}) when $k=3$. 
Similarly, the values for $k=5$ can be obtained from Theorem~\ref{theo:R=1}  
except the three cases where $R=5,6,7$.  The case where $R=6$ follows from 
Theorem~\ref{theo:Rk+1} (or Proposition~\ref{prop:Rk+1}).  
As for the case where $R=5$, $m_5(15Q+5)$ must lie in between $31Q+7$ and $31Q+9$ because 
$m_5(15Q+4)=31Q+6$, $m_5(15Q+6)=31Q+10$ and $m_k(b+1)\geqq m_k(b)+1$ as in Corollary~\ref{coro3}, 
and the value $31Q+8$ would be excluded because $m_k(b)\equiv b\pmod 2$ by Theorem~\ref{mod2}.  
As for the case where $R=7$, the same argument shows that $m_5(15Q+7)=31Q+11, 13$ or $15$.  But 
the value $31Q+15$ would be excluded by \eqref{eq:upper}.  
A similar argument shows the values above when $k=6$.  In fact we also use Proposition~\ref{lowerimp} proved later for 
$R=12, 13, 14$ and $15$.  

Finally we note that $m_5(5)=7$ and $m_6(6)=8$ by Theorem~\ref{mkk0} although we could not determine the values of $m_5(15Q+5)$ 
and $m_6(31Q+6)$ for $Q\geqq 1$ as shown above. 
\end{exam}

\section{Upper and lower bounds of $m_k(b)$} \label{sect:7}

We continue to use the expression
\begin{equation*} \label{QR}
\text{$b=(2^{k-1}-1)Q+R$} 
\end{equation*}
where $Q$ and $R$ are non-negative integers and $0\leqq R\leqq 2^{k-1}-2$. 
Here are naive upper and lower bounds of $m_k(b)$.

\begin{lemm} \label{lemm:uplow} 
$(2^{k}-1)Q+R\leqq m_k(b)\leqq (2^k-1)Q+2R$, i.e. if we denote $m_k(b)=(2^k-1)Q+S$, then 
$R\leqq S\leqq 2R$. 
\end{lemm}

\begin{proof}
We take $a_v=Q+R$ for one $v$ and $a_v=Q$ for all other $v$'s.  
These satisfy \eqref{eq2} and $\sum a_v=(2^k-1)Q+R$, proving the lower bound. 
The upper bound is a restatement of the upper bound in \eqref{eq:upper}. 
\end{proof}

\begin{rema}
It easily follows from Lemma~\ref{lemm:uplow} that 
$\lim_{b\to\infty}m_k(b)/b=(2^k-1)/(2^{k-1}-1)$, 
so $m_k(b)$ is approximately $(2^k-1)b/(2^{k-1}-1)$ when $b$ is large. 
\end{rema}

The bounds in Lemma~\ref{lemm:uplow} are best possible in the sense that 
both $S=R$ and $S=2R$ occur and it is easy to see when $S=R$ occurs.  
In this section we improve the lower bound in Lemma~\ref{lemm:uplow} and see when the lower and upper bounds 
are attained.   
%
%
The following answers the question of when $S=R$ occurs. 

\begin{prop} \label{theo:S=R}
Let $b=(2^{k-1}-1)Q+R$ and $m_k(b)=(2^k-1)Q+S$.  Then $S=R$ if and only if 
$R\leqq k-2$.
\end{prop}

\begin{proof}
The \lq\lq if part" follows from Theorem~\ref{theo:Rk-1}.  
Suppose $R\geqq k-1$.  Then it follows from Lemma~\ref{coro3}, \eqref{eq:divisi} and Theorem~\ref{mkk} that 
\[
\begin{split}
(2^k-1)Q+S&=m_k(b)=m_k((2^{k-1}-1)Q+R)\\
&\geqq m_k(2^{k-1}-Q)+m_k(k-1)+m_k(R-k+1)\\
&\geqq (2^k-1)Q+(k+1)+(R-k+1)\\
&=(2^k-1)Q+R+2
\end{split}
\]
and hence $S\geqq R+2$, proving the \lq\lq only if" part. 
\end{proof}

We shall study when $S=2R$ occurs and improve the lower bound in Lemma~\ref{lemm:uplow} 
in the rest of this section. 
Remember that the polyhedron $P(b)$ defined by $(2^k-1)$ inequalities 
\[
\sum_{(u,v)=0} a_v\leqq b\quad \text{for each $u\in (\Z/2)^k\backslash\{0\}$}
\]
has the point $x=(a_v)$ with $a_v=b/(2^{k-1}-1)$ as the unique vertex and 
the $(2^k-1)$ hyperplanes 
\[
H^u(b)=\{(a_v)\in \R^{2^k-1}\mid \sum_{(u,v)=0} a_v= b\} \quad\text{for $u\in (\Z/2)^k\backslash\{0\}$} 
\]
are in general position.  
We set 
\[
H(m)=\{(a_v)\in\R^{2^k-1}\mid \sum a_v=m\}.
\]
Lemma~\ref{lemm:max} tells us that the intersection 
$P(b)\cap H(m)$ is non-empty if and only if $m\leqq (2^k-1)b/(2^{k-1}-1)$, and that it is 
the one point $x$ if $m=(2^k-1)b/(2^{k-1}-1)$ and a simplex of dimension $2^k-2$ if $m< (2^k-1)b/(2^{k-1}-1)$. 

\begin{lemm} \label{lemm5-1}
Let $u\in (\Z/2)^k\backslash\{0\}$.  Then the $v$-th coordinate $a_v^u$ 
of a vertex $P^u=H(m)\cap(\cap_{u'\not=u}H^{u'})$ of $P(b)\cap H(m)$ is given by 
\[
a_v^u=\begin{cases} 2b-m+(m-b)/2^{k-2} \quad&\text{if $(u,v)\not=0$},\\
m-2b \quad&\text{if $(u,v)=0$}.
\end{cases}
\]
In other words, if $b=(2^{k-1}-1)Q+R$ and $m=(2^k-1)Q+S$, then 
\[
a_v^u=\begin{cases} Q+2R-S+(S-R)/2^{k-2}\quad&\text{if $(u,v)\not=0$},\\
Q+S-2R \quad&\text{if $(u,v)=0$}.
\end{cases}
\]
\end{lemm}

\begin{proof}
Fix $u\in (\Z/2)^k\backslash\{0\}$.  
For each $u'\in (\Z/2)^k\backslash\{0\}$ we consider an equation 
\begin{equation} \label{eq:uv}
\sum_{(u',v')=0}a^u_{v'}=b
\end{equation}
where $v'$ runs over elements with $(u',v')=0$ in the sum.    

The following argument is similar to the latter half of the proof of Lemma~\ref{lemm:max}.  
For $v$ with $(u,v)\not=0$,  we take sum of \eqref{eq:uv} 
over all non-zero $u'$ with $(u',v)=0$.  Then we obtain
\begin{equation} \label{eq:uv1}
(2^{k-1}-1)a^u_v+(2^{k-2}-1)\sum_{v'\not=v}a^u_{v'}= (2^{k-1}-1)b.
\end{equation} 
(Note that $a^u_{v'}$ with $v'\not=v$ appears in the equation \eqref{eq:uv} for $u'$ with 
$(u',v)=(u',v')=0$, so it appears  $(2^{k-2}-1)$ times.)  
Since $a^u_v+\sum_{v'\not=v}a^u_{v'}=m$, we plug $\sum_{v'\not=v}a^u_{v'}=m-a^u_v$ in \eqref{eq:uv1} to obtain 
\begin{equation} \label{eq5-0}
\begin{split}
a^u_v&= \frac{1}{2^{k-2}}\big\{(2^{k-1}-1)b-(2^{k-2}-1)m\big\}\\
&=2b-m+\frac{1}{2^{k-2}}(m-b).
\end{split}
\end{equation}

For $v$ with $(u,v)=0$, we take sum of \eqref{eq:uv} over all non-zero $u'$ with 
$(u',v)=0$ and $u'\not=u$.  Since the number of such $u'$ is $2^{k-1}-2$, we obtain 
\begin{equation} \label{eq5-1}
(2^{k-1}-2)a^u_v+(2^{k-2}-1)\sum_{v'\not=v} a^u_{v'}-\sum_{v'\not=v,(u,v')=0}a^u_{v'}=(2^{k-1}-2)b.
\end{equation}
Here 
\begin{equation} \label{eq5-2}
\sum_{v'\not=v}a^u_{v'}=m-a^u_v 
\end{equation}
and 
\begin{equation} \label{eq5-3}
\begin{split}
\sum_{v'\not=v,(u,v')=0}a^u_{v'}&=m-a^u_v-\sum_{(u,v')\not=0}a^u_{v'}\\
&=m-a^u_v-2^{k-1}\Big(2b-m+\frac{1}{2^{k-2}}(m-b)\Big)\\
&=(2^{k-1}-1)m-(2^{k}-2)b-a^u_v
\end{split}
\end{equation}
where we used \eqref{eq5-0} for $v'$ at the second identity.   
Plugging \eqref{eq5-2} and \eqref{eq5-3} in \eqref{eq5-1}, we obtain 
\[
2^{k-2}a^u_v-2^{k-2}m+(2^k-2)b=(2^{k-1}-2)b
\]
and hence $a^u_v=m-2b$. 
\end{proof}

\begin{prop} \label{prop:upper}
Let $b=(2^{k-1}-1)Q+R$ and $m_k(b)=(2^k-1)Q+S$.  If $S=2R$, then  
$R=2^{k-1}-2^{k-1-\ell}$ for some $0\leqq \ell\leqq k-2$. 
\end{prop}

\begin{proof}
Suppose $S=2R$.  Then it follows from Lemma~\ref{lemm5-1} that the $v$-th coordinate $a_v^u$ of the vertex $P^u$ 
of $P(b)\cap H(m_k(b))$ is given by 
\[
a^u_v=\begin{cases} Q+R/2^{k-2}\quad&\text{if $(u,v)\not=0$},\\
Q \quad&\text{if $(u,v)=0$}. 
\end{cases}
\]
Since $m_k(b)=(2^k-1)Q+S$ and $S=2R$ by assumption, 
there is a lattice point on the simplex $P(b)\cap H(m_k(b))$.  The simplex is the convex hull of the vertices $P^u$, so there exist 
non-negative real numbers $t_u$'s with $\sum t_u=1$ such that $\sum t_uP^u$ is a lattice point, i.e.  
\[
\sum t_ua^u_v=\sum_{(u,v)\not=0}t_u(Q+R/2^{k-2})+\sum_{(u,v)=0}t_uQ=Q+(\sum_{(u,v)\not=0}t_u)R/2^{k-2}\in \Z
\]
for any $v$. This means that $(\sum_{(u,v)\not=0}t_u)R/2^{k-2}=0$ or $1$, i.e.  
\begin{equation} \label{eq:5-1}
\sum_{(u,v)\not=0}t_u=0\quad\text{or}\quad 2^{k-2}/R\quad\text{for any $v$}
\end{equation}
because $0\leqq R\leqq 2^{k-1}-2$ and $\sum_{(u,v)\not=0}t_u\leqq 1$. 
On the other hand, 
\begin{equation} \label{eq:5-2}
\sum_v\sum_{(u,v)\not=0}t_u=2^{k-1}
\end{equation}
because each $t_u$ appears $2^{k-1}$ times in the sum above and $\sum t_u=1$.  
It follows from \eqref{eq:5-1} and \eqref{eq:5-2} that there are exactly $2R$ numbers of $v$'s such that 
$\sum_{(u,v)\not=0}t_u\not=0$, in other words, there are exactly $2^k-1-2R$ numbers of 
$v$'s such that $\sum_{(u,v)\not=0}t_u=0$.  The identity $\sum_{(u,v)\not=0}t_u=0$ implies that 
$t_u=0$ for all $u$ with $(u,v)\not=0$ since $t_u\geqq 0$.  
Based on these observations, we introduce 
\[
\begin{split}
U:=&\text{ the linear span of }U_0:=\{u\mid t_u\not=0\},\\
V:=&\text{ the linear span of }V_0:=\{v \mid t_u=0 \quad\text{for $\forall u$ such that $(u,v)\not=0$}\}. 
\end{split}
\]
If $v\in V_0$, then it follows from the definition of $U_0$ and $V_0$ that $(u,v)=0$ for any $u\in U_0$ and hence 
$(u,v)=0$ for any $u\in U$ since $U$ is the linear span of $U_0$.  This implies that $(u,v)=0$ for any $u\in U$ and 
$v\in V$ since $V$ is the linear span of $V_0$.  It follows that 
\begin{equation} \label{eq:5-3}
\dim U\leqq k-\dim V.
\end{equation}
We note that $V$ contains at least $2^k-1-2R$ non-zero elements by the observation made above. 

Suppose that 
\begin{equation} \label{eq:5-4}
\text{$2^{k-1}-2^{k-1-\ell} \leqq R<  2^{k-1}-2^{k-1-(\ell+1)}$ for some $0\leqq\ell \leqq k-2$.}
\end{equation}
(Note that $R$ lies in the inequality \eqref{eq:5-4} for some $\ell$ because $0\leqq R\leqq 2^{k-1}-2$.) 
Then, since 
$2^{k-\ell-1}-1< 2^k-1-2R$ and $V$ contains at least $2^k-1-2R$ non-zero elements, 
$V$ contains at least $2^{k-\ell-1}$ non-zero elements and hence $\dim V\geqq k-\ell$.  This together with \eqref{eq:5-3} 
shows  
\begin{equation} \label{eq:5-4-1}
\dim U\leqq\ell
\end{equation} 

Since the bilinear form $(\ ,\ )$ is non-degenerate, there is a subspace $W$ of $(\Z/2)^k$ such that $\dim W=\dim U$ 
and the bilinear form $(\ ,\ )$ restricted to $U\times W$ is still non-degenerate.  
We take sum of \eqref{eq:5-1} over all non-zero $v\in W$. 
In this sum, each $t_u$ for $u\in U\backslash\{0\}$ appears $2^{\dim W-1}$ times. 
Since $\dim W=\dim U$ and $\sum_{u\in U\backslash\{0\}}t_u=1$, we obtain 
\[
2^{\dim U-1}\leqq (2^{\dim U}-1)2^{k-2}/R
\]
and hence 
\begin{equation} \label{eq:5-5}
R\leqq (2^{\dim U}-1)2^{k-\dim U-1}\leqq 2^{k-1}-2^{k-\ell-1}
\end{equation}
where we used \eqref{eq:5-4-1} at the latter inequality.  
Then \eqref{eq:5-4} and \eqref{eq:5-5} show that $R=2^{k-1}-2^{k-1-\ell}$, proving 
the proposition.  
\end{proof}

It turns out that the converse of Proposition~\ref{prop:upper} holds, i.e. $S=2R$ can be attained when $R=2^{k-1}-2^{k-1-\ell}$. 
In fact, we can prove the following. 

\begin{prop} \label{prop:lower}
Let $b=(2^{k-1}-1)Q+R$ and let $2^{k-1}-2^{k-1-\ell}\leqq R<2^{k-1}-2^{k-1-(\ell+1)}$ 
for some $0\leqq \ell \leqq k-2$. 
Then 
$$m_k(b)\geqq (2^k-1)Q+R+2^{k-1}-2^{k-1-\ell}.$$
In particular, if $R=2^{k-1}-2^{k-1-\ell}$ for some $0\leqq \ell\leqq k-2$, then 
$m_k(b)\geqq (2^k-1)Q+2R$.
\end{prop}

\begin{proof}
We take 
$$m=(2^k-1)Q+R+2^{k-1}-2^{k-1-\ell}$$ 
and find a lattice point in the simplex $P(b)\cap H(m)$ with non-negative coordinates. 
Set 
$$r=R-2^{k-1}+2^{k-1-\ell}.$$
The $v$-th coordinate $a^u_v$ of the vertex $P^u$ of $P(b)\cap H(m)$ is given by  
\begin{equation} \label{eq:5-6}
a^u_v=\begin{cases} Q+r+2-2^{1-\ell} \quad&\text{if $(u,v)\not=0$},\\
Q-r \quad&\text{if $(u,v)=0$} \end{cases}
\end{equation}
by Lemma~\ref{lemm5-1}.  Set 
\begin{equation} \label{eq:5-7}
L=2-2^{1-\ell}.
\end{equation}
Any point in $P(b)\cap H(m)$ can be expressed as $\sum_{u\in (\Z/2)^k\backslash\{0\}} t_uP^u$ 
with $t_u\geqq 0$ and $\sum t_u=1$, and we find from \eqref{eq:5-6} that its $v$-th coordinate $a_v$ is given by 
\begin{equation} \label{eq:5-8}
\begin{split}
a_v=&(\sum_{(u,v)\not=0}t_u)(Q+r+L)+(\sum_{(u,v)=0}t_u)(Q-r)\\
=&(\sum t_u)Q+(1-\sum_{(u,v)=0}t_u)(r+L)+(\sum_{(u,v)=0}t_u)(-r)\\
=&Q+r+L-(\sum_{(u,v)=0}t_u)(2r+L).
\end{split}
\end{equation}

We take a codimension 1 subspace $V$ of $(\Z/2)^k$ and an $\ell$-dimensional subspace $U$ of $V$ arbitrarily and 
define 
\begin{equation} \label{eq:5-8-1}
t_u=\begin{cases} \frac{2r}{2r+L}\frac{1}{2^{k-1}} \quad&\text{for $u\notin V$},\\
\frac{L}{2r+L}\frac{1}{2^\ell-1} \quad&\text{for $u\in U\backslash\{0\}$},\\
0 \quad&\text{otherwise}.\end{cases}
\end{equation}
Then $t_u\geqq 0$ and $\sum t_u=1$.  
We shall check that $a_v$ in \eqref{eq:5-8} is a non-negative integer.  
We denote by $v^\perp$ the codimension 1 subspace of $(\Z/2)^k$ consisting of elements $w$ such that $(v,w)=0$ 
and distinguish three cases according to the position of $v^\perp$ relative to $V$ and $U$.  

{\bf Case 1.} The case where $v^\perp=V$.  In this case, 
\[
\sum_{(u,v)=0}t_u=\frac{L}{2r+L}\frac{1}{2^{\ell}-1}(2^\ell-1)=\frac{L}{2r+L},
\]
so $a_v=Q+r$ by \eqref{eq:5-8}. 

{\bf Case 2.} The case where $v^\perp\not=V$ and $v^\perp\supset U$.  In this case, 
$v^\perp\cap V$ is of dimension $k-2$ and 
\[
\sum_{(u,v)=0}t_u=\frac{2r}{2r+L}\frac{1}{2^{k-1}} 2^{k-2}+\frac{L}{2r+L}\frac{1}{2^{\ell}-1}(2^\ell-1)=\frac{r+L}{2r+L},
\]
so $a_v=Q$ by \eqref{eq:5-8}. 

{\bf Case 3.} The case where $v^\perp\not=V$ and $v^\perp\not\supset U$. 
In this case, $v^\perp\cap V$ is of dimension $k-2$ and $v^\perp\cap U$ is of dimension $\ell-1$ and hence 
\[
\begin{split}
\sum_{(u,v)=0}t_u&=\frac{2r}{2r+L}\frac{1}{2^{k-1}}2^{k-2}+\frac{L}{2r+L}\frac{1}{2^\ell-1}(2^{\ell-1}-1)\\
&=\frac{r+L-1}{2r+L}
\end{split}
\]
where we used \eqref{eq:5-7} at the second identity, so $a_v=Q+1$ by \eqref{eq:5-8}.  

In any case $a_v$ is a non-negative integer, so $\sum_{u\in (\Z/2)^k\backslash\{0\}} t_uP^u$ 
with $t_u$ in \eqref{eq:5-8-1} is a lattice point in $P(b)\cap H(m)$ with non-negative coordinates.  
This proves the proposition.
\end{proof}

Now we are ready to prove the latter theorem in the Introduction. 

\begin{theo} \label{theo:main}
Let $b=(2^{k-1}-1)Q+R$.
If $2^{k-1}-2^{k-1-\ell}\leqq R<2^{k-1}-2^{k-1-(\ell+1)}$ for some $0\leqq \ell\leqq k-2$, then 
\[
(2^k-1)Q+R+2^{k-1}-2^{k-1-\ell} \leqq m_k(b) \leqq  (2^k-1)Q+2R
\]
where the lower bound is attained if and only if $R-(2^{k-1}-2^{k-1-\ell})\leqq k-\ell-2$ and 
the upper bound is attained if and only if $R=2^{k-1}-2^{k-1-\ell}$.  
\end{theo}

\begin{proof}
The inequality and the statement on the upper bound follows from Propositions~\ref{prop:upper} and~\ref{prop:lower}. 
Moreover, Theorem~\ref{theo:R=1} shows that the lower bound is attained if $R-(2^{k-1}-2^{k-1-\ell})\leqq k-\ell-2$. 
Suppose $R-(2^{k-1}-2^{k-1-\ell})\geqq k-\ell-1$ and set 
\begin{equation} \label{D}
D=R-(2^{k-1}-2^{k-1-\ell})-(k-\ell-1).  
\end{equation}
Then it follows from Lemma~\ref{coro3} and Theorem~\ref{theo:R=1} that 
\[
\begin{split}
m_k(b)&=m_k((2^{k-1}-1)Q+R)\\
&=m_k((2^{k-1}-1)Q+2^{k-1}-2^{k-1-\ell}+k-\ell-1+D)\\
&\geqq m_k((2^{k-1}-1)Q+2^{k-1}-2^{k-1-\ell}+k-\ell-1)+m_k(D)\\
&\geqq (2^k-1)Q+2^k-2^{k-\ell}+k-\ell+1+D\\
&=(2^k-1)Q+R+2^{k-1}-2^{k-\ell-1}+2
\end{split}
\]
where we used \eqref{D} at the last identity.   Therefore the lower bound is not attained if $R-(2^{k-1}-2^{k-1-\ell})\geqq k-\ell-1$. 
\end{proof}

\section{A slight improvement of lower bounds}

When $R\leqq 2^{k-2}-1$, the lower bound of $m_k(b)$ in Theorem~\ref{theo:main} is nothing but $(2^k-1)Q+R$ and this is 
an obvious lower bound.  In this section we improve the lower bound when $2^{k-2}-4\leqq R\leqq 2^{k-2}-1$.  

\begin{prop} \label{lowerimp}
If $k$ is odd, then 
\begin{enumerate}
\item $m_k(2^{k-1}-1)Q+2^{k-2}-1)\geqq  (2^k-1)Q+2^{k-1}-k$,
\item $m_k(2^{k-1}-1)Q+2^{k-2}-2)\geqq (2^k-1)Q+2^{k-1}-k-1$.
\end{enumerate}
If $k$ is even, then 
\begin{enumerate}
\item $m_k(2^{k-1}-1)Q+2^{k-2}-1)\geqq  (2^k-1)Q+2^{k-1}-k+1$,
\item $m_k(2^{k-1}-1)Q+2^{k-2}-2)\geqq (2^k-1)Q+2^{k-1}-k-2$,
\item $m_k(2^{k-1}-1)Q+2^{k-2}-3)\geqq (2^k-1)Q+2^{k-1}-2k+1$,
\item $m_k(2^{k-1}-1)Q+2^{k-2}-4)\geqq (2^k-1)Q+2^{k-1}-2k$. 
\end{enumerate}
\end{prop}

\begin{proof}
In any case it suffices to prove the inequality when $Q=0$ by Lemma~\ref{coro3} (2). 
We recall how $m_k(2^{k-2})=2^{k-1}$ is obtained.  
Choose any non-zero element $u_0\in (\Z/2)^k$ and define 
\begin{equation} \label{avv}
a_v=\begin{cases} 1 \quad&\text{if $(u_0,v)\not=0$,}\\
0\quad&\text{if $(u_0,v)=0$.} \end{cases}
\end{equation}
Then 
\[
\sum_{(u,v)=0}a_v=\begin{cases} 2^{k-2} \quad&\text{if $u\not=u_0$,}\\
0\quad&\text{if $u=u_0$}\end{cases}
\]
and $\sum a_v=2^{k-1}$.  This attains $m_k(2^{k-2})=2^{k-1}$.  

We take $$u_0=(1,\dots,1)^t.$$  Then $(u_0,v)=0$ if and only if the number of 1 in the components of $v$ 
is even.  Let 
\[
\begin{split}
V_1:&=\{\e_1,\dots,\e_k\} \subset (\Z/2)^k\\
V_2:&=\begin{cases} 
V_1\cup\{ u_0\}\quad&\text{for $k$ odd},\\
V_1\cup\{ u_0-\e_1, u_0-\e_2\}\quad&\text{for $k$ even},\end{cases}
\end{split}
\]
and define for $q=1,2$ 
\[
a_v^{(q)}:=\begin{cases} 1 \quad&\text{if $(u_0,v)\not=0$ and $v\notin V_q$,}\\
0 \quad&\text{otherwise.}\end{cases}
\]
One can check that $\sum_{(u,v)=0}a_v^{(q)}\leqq 2^{k-2}-q$ for any non-zero $u\in (\Z/2)^k$.  
Clearly  
\[
\sum a_v^{(q)}=\begin{cases} 2^{k-1}-k \quad&\text{when $q=1$},\\
2^{k-1}-k-1 \quad&\text{when $q=2$ and $k$ is odd},\\
2^{k-1}-k-2 \quad&\text{when $q=2$ and $k$ is even}. \end{cases}
\]  
This together with the congruence $m_k(b)\equiv b\pmod 2$ in Theorem~\ref{mod2} (applied when $q=1$ and $k$ is even) 
implies the inequalities (1) and (2) 
in the proposition.   

The proof of the inequality (4) is similar.  Assume $k$ is even and let 
\[
V_4:=V_1\cup\{ u_0-\e_1, \dots, u_0-\e_k\}
\]
and define 
\[
a_v^{(4)}:=\begin{cases} 1 \quad&\text{if $(u_0,v)\not=0$ and $v\notin V_4$,}\\
0 \quad&\text{otherwise.}\end{cases}
\]
One can check that $\sum_{(u,v)=0}a_v^{(4)}\leqq 2^{k-2}-4$ for any 
non-zero $u\in (\Z/2)^k$ (where we use the assumption on $k$ being even) and $\sum a_v^{(4)}=2^{k-1}-2k$.  Therefore 
\begin{equation*} 
\text{$m_k(2^{k-2}-4)\geqq 2^{k-1}-2k$}  
\end{equation*}
which implies the inequality (4) in the proposition.  The inequality (3) follows from (4) since $m_k(b+1)\geqq 
m_k(b)+1$.  
\end{proof}

\section{Some observation on the Conjecture} \label{sect:8}

The Conjecture in Section~\ref{sect:5} says that 
\[
m_k((2^{k-1}-1)Q+R)=(2^k-1)Q+m_k(R)
\]
and this is equivalent to saying
\begin{equation} \label{conj}
m_k(b+2^{k-1}-1)=m_k(b)+2^k-1.
\end{equation}
In this section, we prove \eqref{conj} when $b$ is large, to be more precise,  we prove the following.

\begin{theo} \label{theo5-1}
Let $b=(2^{k-1}-1)Q+R$.  If 
\[
Q\geqq \begin{cases} R \quad&\text{when\quad $0\leqq R\leqq 2^{k-2}-1$},\\
R-2^{k-2} \quad&\text{when\quad $2^{k-2}\leqq R\leqq 2^{k-1}-2$},
\end{cases}
\]
{\rm (}this is the case when 
$b\geqq  (2^{k-1}-1)(2^{k-2}-1)${\rm )}, then 
\[
m_k(b+2^{k-1}-1)=m_k(b)+2^k-1.
\]
\end{theo}

\begin{proof}
By Lemma~\ref{coro3} (2), it suffices to prove 
\begin{equation} \label{eq5-4}
m_k(b+2^{k-1}-1)\leqq m_k(b)+2^k-1.
\end{equation}
Remember the polyhedron $P(b)$ defined by $(2^k-1)$ inequalities 
\begin{equation} \label{eq:6-1}
\sum_{(u,v)=0} a_v\leqq b\quad \text{for each $u\in (\Z/2)^k\backslash\{0\}$}.
\end{equation}
We will find $m$ such that the intersection of $P(b+2^{k-1}-1)$ with a half space $H^+(m)$ in $\R^{2^k-1}$ defined by 
\[
H^+(m)=\{ \sum a_v\geqq m\}
\] 
has a lattice point with coordinates $\geqq 1$. 

{\bf Case 1.} The case where $0\leqq R\leqq 2^{k-2}-1$.  In this case we take $$m=(2^k-1)(Q+1)+R.$$
Since $$b+2^{k-1}-1=(2^{k-1}-1)(Q+1)+R,$$
the coordinates of a vertex (except the vertex $x$ of $P(b+2^{k-1}-1)$) in $P(b+2^{k-1}-1)\cap H^+(m)$ are either 
$Q+1+R$ or $Q+1-R$ by Lemma~\ref{lemm5-1}, so those vertices are lattice points and their coordinates are 
greater than or equal to $1$ since $Q\geqq R$ by assumption.  We know  
$$m_k(b+2^{k-1}-1)\geqq (2^k-1)(Q+1)+R$$  
by Lemma~\ref{lemm:uplow}, so any lattice point 
$(a_v)$ in \eqref{eq:6-1} with $b$ replaced by $b+2^{k-1}-1$, at which $\sum a_v$ attains the maximum value $m_k(b+2^{k-1}-1)$,  
lies in $P(b+2^{k-1}-1)\cap H^+(m)$ and hence $a_v\geqq 1$ for every $v$.  Since   
$\{a_v-1\}$ is a set of non-negative integers which satisfy \eqref{eq:6-1} and $$\sum(a_v-1)=m_k(b+2^{k-1}-1)-(2^k-1),$$
it follows from the definition of $m_k(b)$ that 
$$m_k(b+2^{k-1}-1)-(2^k-1)\leqq m_k(b),$$ proving the desired inequality \eqref{eq5-4}.  

{\bf Case 2.} The case where $2^{k-2}\leqq R\leqq 2^{k-1}-2$.  In this case we take 
\[
m=(2^k-1)(Q+1)+R+2^{k-2}.
\]
Then the coordinates of a vertex (except the vertex $x$) in $P(b+2^{k-1}-1)\cap H^+(m)$ are either 
$Q+2+R-2^{k-2}$ or $Q+1-R+2^{k-2}$ by Lemma~\ref{lemm5-1}, so those vertices are lattice points and their coordinates are 
greater than or equal to $1$ 
since $Q\geqq R-2^{k-2}$ by assumption.  We know 
$$m_k(b+2^{k-1}-1)\geqq (2^k-1)(Q+1)+R+2^{k-2}$$  
by Proposition~\ref{prop:lower}, so any lattice point 
$(a_v)$ in \eqref{eq:6-1} with $b$ replaced by $b+2^{k-1}-1$, at which $\sum a_v$ attains the maximum value $m_k(b+2^{k-1}-1)$,  
lies in $P(b+2^{k-1}-1)\cap H^+(m)$ and hence $a_v\geqq 1$ for every $v$. The remaining argument is same as in Case 1 above.  
\end{proof}

\newpage

\section*{Appendix}

Below is a table of values of $s_\R(m,p)$ for $2\leqq p\leqq 18$ and $2\leqq m\leqq 40$.  

{\Tiny
\begin{table}[htb]
\begin{tabular}{|c||c|c|c|c|c|c|c|c|c|c|c|c|c|c|c|c|c|} \hline
$m\backslash p$ & 2 & 3 & 4 & 5 & 6 & 7 & 8 & 9 & 10 &11 & 12 & 13 & 14 & 15 & 16 & 17 & 18   \\
\hline\hline 2 & 2 &    &    &     &   &      &     &    &    &     &     &    &     &     &     &     &       \\
\hline  3         & 2 & 3 &    &    &    &      &     &    &    &     &     &    &     &     &     &     &       \\
\hline 4          & 1&  3 & 4 &    &    &      &     &    &    &     &     &    &     &     &     &     &       \\
\hline 5          & 1 & 2 & 4 & 5 &    &      &     &    &    &     &     &    &     &     &     &     &        \\
\hline 6          & 1 & 2 & 3 & 5 & 6  &     &     &    &    &     &     &    &     &     &     &     &        \\
\hline 7          & 1 & 1 & 3 & 4 & 6  & 7  &     &    &    &     &     &    &     &     &     &     &        \\
\hline 8          & 1 & 1 & 2 & 4 & 4  & 7  & 8  &    &    &     &     &    &     &     &     &     &        \\
\hline 9          & 1 & 1 & 2 & 2 & 4  & 5  & 8  &  9 &    &     &    &     &    &     &     &     &        \\
\hline 10        & 1 & 1 & 1 & 2 & 3 &  5  &  6 &  9 & 10 &   &    &     &     &    &      &    &        \\
\hline 11        & 1 & 1 & 1 & 2 & 3  & 4  &  6 &  7 & 10 &11&   &     &     &     &      &    &      \\
\hline 12        & 1 & 1 & 1 & 2 & 2  & 4  & $*\leqq 5$  & 7  & 8  &11& 12&     &     &     &      &     &  \\
\hline 13        & 1 & 1 & 1 & 1 & 2  & 3  & $*$  & $*\leqq 6$  & 8  & 9 & 12& 13&    &      &          &    &  \\
\hline 14        & 1 & 1 & 1 & 1 & 2  & 3  & 4  &  $*$ & $*\leqq 7$  & 9 & 10& 13& 14 &    &          &     &   \\
\hline 15        & 1 & 1 & 1 & 1 & 2  & 2  & 4  & 5 &  *  & $*\leqq 8$ & 10& 11& 14 & 15 &    &      &   \\
\hline 16        & 1 & 1 & 1 & 1 & 1  & 2  & 2  & 5  & *  & * & $*\leqq 9$ & 11& 11 & 15 & 16 &     &   \\
\hline 17        & 1 & 1 & 1 & 1 & 1  & 2  & 2  & 3  & $*\geqq 5$  & * &  * & $*\leqq 10$  & 11 & 12 & 16 & 17&   \\
\hline 18        & 1 & 1 & 1 & 1 & 1  & 2  & 2  & 3  & 3  &$*\geqq 5$  &  * &  * &  * &  12&  13& 17 & 18   \\
\hline 19        & 1 & 1 & 1 & 1 & 1  & 1  & 2  & 2  & 3  & 4 &  * &  *  & *  & *  & 13 & 14 & 18  \\
\hline 20        & 1 & 1 & 1 & 1 & 1  & 1  & 2  & 2  & 3  & 4 &  5 &  *  & *  & *  &  * &  14&  15  \\
\hline 21        & 1 & 1 & 1 & 1 & 1  & 1  & 2  & 2  & 3  & 3 &  5 &  *  & *  & *  &  * &  *  &   15  \\
\hline 22        & 1 & 1 & 1 & 1 & 1  & 1  & 1  & 2  & 2  & 3 &  4 &  *  & *  & *  &  * &  *  &  *      \\
\hline 23        & 1 & 1 & 1 & 1 & 1  & 1  & 1  & 2  & 2  & 2 &  4 &  5  & *  & *  &  * &  *  &  *        \\
\hline 24        & 1 & 1 & 1 & 1 & 1  & 1  & 1  & 2  & 2  & 2 &  3 &  5  & *  & *  &  *    & *    & *  \\
\hline 25        & 1 & 1 & 1 & 1 & 1  & 1  & 1  & 1  & 2  & 2 &  3 &  3  & $*\geqq 5$  & *  &  *    &   *  & *  \\
\hline 26        & 1 & 1 & 1 & 1 & 1  & 1  & 1  & 1  & 2  & 2 &  2 &  3  & 4 &  * &  * &  *  & *         \\
\hline 27        & 1 & 1 & 1 & 1 & 1  & 1  & 1  & 1  & 2  & 2 &  2 &  3  & 4  & *  &  * & *   & *        \\
\hline 28        & 1 & 1 & 1 & 1 & 1  & 1  & 1  & 1  & 1  & 2 &  2 &  3  & 3  & $*\geqq 5$  &  *    &  *     & *  \\
\hline 29        & 1 & 1 & 1 & 1 & 1  & 1  & 1  & 1  & 1  & 2 &  2 &  2  & 3  & 4  &  * &  *  & *        \\
\hline 30        & 1 & 1 & 1 & 1 & 1  & 1  & 1  & 1  & 1  & 2 &  2 &  2  & 2  & 4  &  5 & *   & *        \\
\hline 31        & 1 & 1 & 1 & 1 & 1  & 1  & 1  & 1  & 1  & 1 &  2 &  2  & 2  & 3  &  5 &   6    & *  \\
\hline 32        & 1 & 1 & 1 & 1 & 1  & 1  & 1  & 1  & 1  & 1 &  2 &  2  & 2  & 3  &  3 &   6    &  *       \\
\hline 33        & 1 & 1 & 1 & 1 & 1  & 1  & 1  & 1  & 1  & 1 &  2 &  2  & 2  & 2  &  3 & 3 & $*\geqq 6$        \\
\hline 34        & 1 & 1 & 1 & 1 & 1  & 1  & 1  & 1  & 1  & 1 &  1 &  2  & 2  & 2  &  3 & 3 &   4   \\
\hline 35        & 1 & 1 & 1 & 1 & 1  & 1  & 1  & 1  & 1  & 1 &  1 &  2  & 2  & 2  &  3 & 3 &   4  \\
\hline 36        & 1 & 1 & 1 & 1 & 1  & 1  & 1  & 1  & 1  & 1 &  1 &  2  & 2  & 2  &  2 & 3 &   3   \\
\hline 37        & 1 & 1 & 1 & 1 & 1  & 1  & 1  & 1  & 1  & 1 &  1 &  1  & 2  & 2  &  2 & 2 &   3  \\
\hline 38        & 1 & 1 & 1 & 1 & 1  & 1  & 1  & 1  & 1  & 1 &  1 &  1  & 2  & 2  &  2 & 2 &   3  \\
\hline 39        & 1 & 1 & 1 & 1 & 1  & 1  & 1  & 1  & 1  & 1 &  1 &  1  & 2  & 2  &  2 & 2 &   3  \\
\hline 40        & 1 & 1 & 1 & 1 & 1  & 1  & 1  & 1  & 1  & 1 &  1 &  1  & 1  & 2  &  2 & 2 &   2  \\
\hline
\end{tabular} \vspace{1.25ex}
\caption{$s_\R(m,p)$ for $2\leqq p\leqq 18,\ 2\leqq m\leqq 40$}
\end{table}
}
Since $s_\R(m,1)=1$, the case where $p=1$ is omitted.  
Remember that $s_\R(m,p)=1$ if and only if $m\geqq 3p-2$ by Theorem~\ref{theo3} and that 
the values of $s_\R(m,p)$ for $p=m-1,m-2$ and $m-3$ can be obtained from Theorems~\ref{theo1} and~\ref{theo2}.  
The other values can be obtained from Table~\ref{table1} in Section~\ref{sect:56} and the fact that 
$s_\R(m,p)=k$ for $k\geqq 2$ if and only if $m_{k+1}(p-1)<m\leqq m_k(p-1)$ (Lemma~\ref{lemm:sR}).  
The asterisk $*$ in a box means that the value is unknown. 
Finally we note that $s_\R(m,p)$ increases as $p$ increases while it decreases as $m$ increases (Proposition~\ref{prop2}).


\begin{thebibliography}{19}


\bibitem{aize09}
A.~Aizenberg, Graduate thesis, Moscow State University 2009.

\bibitem{buch08}
V.~M.~Buchstaber,
\emph{Lectures on toric topology}, 
Trends in Mathematics, Information Center for Mathematical Sciences, 
vol. 11, No. 1 (2008), 1--55. 

\bibitem{bu-pa02}
V.~M.~Buchstaber and T.~E.~Panov,
\emph{Torus Actions and Their Applications in Topology and Combinatorics},
University Lecture, vol.~{\bf24},
Amer. Math. Soc., Providence, R.I., 2002.

\bibitem{da-ja91}
M. W. Davis and T. Januszkiewicz,
\emph{Convex polytopes, Coxeter orbifolds and torus actions},
Duke Math. J.~{\bf62} (1991), 417--451.

\bibitem{erok09}
N. Erokhovets, 
\emph{Buchstaber invariant of simple polytopes}, 
arXiv:0908.3407. 

\bibitem{cont08}
M. Harada, Y. Karshon, M. Masuda and T. Panov (eds),
\emph{Toric Topology}, Proc. of  the International Conference held at Osaka City University in 2006,  
Contemp. Math. 460 (2008).  

\end{thebibliography}
\end{document}